\documentclass[11pt,a4paper]{article}
\usepackage{amssymb,amsfonts,amsmath,amsthm}
\usepackage[latin1]{inputenc}
\usepackage[english]{babel}
\usepackage[left=2.7cm,right=2.7cm,bottom=2.67cm,top=2.67cm]{geometry}
\usepackage{dsfont}
\usepackage{mathrsfs}
\usepackage{natbib}
\usepackage{enumerate}
\usepackage{color}
\usepackage{sectsty}
\RequirePackage[OT1]{fontenc}
\usepackage{graphicx}
\usepackage{multirow}
\usepackage{multicol}
\usepackage{subfigure}
\usepackage{placeins}
\usepackage{rotating}
\usepackage[flushleft]{threeparttable}
\usepackage{booktabs}
\usepackage[normalem]{ulem}

\graphicspath{{fig/}}

\theoremstyle{plain}

\newtheorem{lema}{Lemma}

\theoremstyle{definition}

\newcommand{\E}{\mathds{E}}
\newcommand{\F}{\mathscr{F}}
\newcommand{\R}{\mathds{R}}

\newcommand{\Z}{\mathds{Z}}

\newcommand{\bs}[1]{\boldsymbol{#1}}
\renewcommand{\qed}{\hfill \mbox{\raggedright \rule{.07in}{.1in}}}

\long\def\sfootnote[#1]#2{\begingroup%
\def\thefootnote{\fnsymbol{footnote}}\footnote[#1]{#2}\endgroup}

\def\bfootnote{\xdef\@thefnmark{}\@footnotetext}

\begin{document}
\pagestyle{myheadings} % allow the use of headings
\markboth{UWARMA}{G. Pumi and T.S. Prass}

\thispagestyle{empty}
{\centering
\Large{\bf Unit-Weibull Autoregressive Moving Average Models}\vspace{.5cm}\\
\normalsize{ {\bf Guilherme Pumi${}^{\mathrm{a,}}$\sfootnote[1]{Corresponding author. This Version: \today},\let\thefootnote\relax\footnote{\hskip-.3cm$\phantom{s}^\mathrm{a}$Mathematics and Statistics Institute and Programa de P\'os-Gradua\c c\~ao em Estat\'istica - Universidade Federal do Rio Grande do Sul.
%- 9500,  Bento Gon\c calves Avenue - 91509-900, Porto Alegre - RS - Brazil.
} Taiane Schaedler Prass${}^\mathrm{a}$ and Cleiton Guollo Taufemback${}^\mathrm{a}$
 \\
\let\thefootnote\relax\footnote{E-mails: guilherme.pumi@ufrgs.br (G. Pumi), taiane.prass@ufrgs.br (T.S. Prass)}
\let\thefootnote\relax\footnote{ORCIDs: 0000-0002-6256-3170 (Pumi); 0000-0003-3136-909X (Prass); 0000-0001-5714-1160 (Taufemback).}}\\
\vskip.3cm
}}
% \affiliation{Universidade Federal do Rio Grande do Sul}
% \address[a]{Department of Statistics and Graduate Program in Statistics \\
% Universidade Federal do Rio Grande do Sul\\
% 91509-900, Porto Alegre\\
% Brazil\\
% \printead{e1,e2}}

\begin{abstract}
In this work we introduce the class of unit-Weibull Autoregressive Moving Average models for continuous random variables taking values in $(0,1)$.
The proposed model is an observation driven one, for which, conditionally on a set of covariates and the process' history, the random component is assumed to follow a unit-Weibull distribution parameterized through its $\rho$th quantile. The systematic component prescribes an ARMA-like structure to model the conditional $\rho$th quantile by means of a link. Parameter estimation in the proposed model is performed using partial maximum likelihood, for which we provide closed formulas for the score vector and partial information matrix. We also discuss some inferential tools, such as the construction of confidence intervals, hypotheses testing, model selection, and forecasting. A Monte Carlo simulation study is conducted to assess the finite sample performance of the proposed partial maximum likelihood approach.  Finally, we examine the prediction power by contrasting our method with others in the literature using the Manufacturing Capacity Utilization from the US.
\vspace{.2cm}

\vspace{.2cm}

\end{abstract}
% Obrigatorio na TEST
\textbf{Statements and Declarations:} The authors declare that they have NO affiliations with or involvement in any
organization or entity with any financial interests in the subject matter or materials discussed in this manuscript.

% . The constructive comments of two Reviewers and the Associate Editor are gratefully acknowledged.

\section{Introduction}

In recent years, the interest in non-Gaussian time series modeling has grown considerably. This is most noticeable for time series supported on the open unitary interval $(0,1)$, as, for instance, relative humidity, the incidence of COVID19 in some particular region (per 10,000 inhabitants, say), the level of a reservoir, etc. Since models supported on the real line may provide forecasts outside the range of plausible values for the phenomenon under study, as in \cite{grande},  specific models are needed to accommodate the natural bounds in the data. This problem typically hinders the application of classical models such as the class of autoregressive moving average (ARMA) and its extensions in the context of double-bounded time series. One approach capable of handling double-bounded time series that has gained attention in the literature over the last decade is the so-called GARMA (generalized ARMA) discussed in \cite{benjamin2003}. In few words, GARMA modeling merges the strengths of ARMA modeling into a generalized linear model framework, yielding a class of very flexible models that can be easily tailored to accommodate a wide variety of structures including non-gaussianity, bounds, asymmetries, etc.

GARMA are classified as observation driven \cite{cox1981} models and, as such,  they are defined by two components: the random component, which specifies the probability structure (conditional distributions) for the model; and the systematic component, which is responsible for the dependence structure present in the model. In practice, original GARMA models only considered random components from the canonical exponential family and an ARMA-like systematic component to model the conditional mean. However, as the literature on GARMA models for continuous double-bounded time series has grown over the years, its scope has grown as well. The interest nowadays lies in models for which the systematic component follows the usual approach of GLM with an additional dynamic term of the form
\begin{equation*}
g(\mu_t)={\eta}_t=\bs{x}_{t}'\bs{\beta}+\tau_t,
\end{equation*}
where $g$ is a suitable link function, $\mu_t$ is some quantity of interest, $\bs{x}_t$ denotes a vector of {(possibly random and time dependent)} covariates observed at time $t$ with associated vector of coefficient $\bs\beta$  and $\tau_t$ is a term responsible to accommodate any serial dependence in $\mu_t$. The terms $\mu_t$ and $\tau_t$ vary according to the model's scope and intended application. For instance, the classical ARMA form and its variants are used in \cite{Rocha2017, Maior,Prass} and \cite{Bayers, Bayerseas} while \cite{Pumi2017} and \cite{helen} apply a long range dependent ARFIMA specification. More exotic non-linear specifications can also be found, as in the case of the Beta autoregressive chaotic models of \cite{BARC}.

Another important feature in GARMA models is the nature of $\mu_t$. Originally, $\mu_t$ denoted the model's conditional mean at time $t$, as in Beta-based models such as \cite{Rocha2017,Bayerseas, Pumi2017} and \cite{BARC} or as in models for positive time series \citep{Prass}. However, for KARMA models, $\mu_t$ denotes the model's median while for SYMARFIMA models, when the random component does not possess a finite first moment, $\mu_t$ denotes the model's point of symmetry. In this work we introduce a new GARMA-like model for which the random component follows a carefully parameterized Unit-Weibull distribution while the systematic component follows an ARMA like structure, resembling  the $\beta$ARMA and KARMA models, but with an important difference: it models the model's $\rho$th quantile, for any $\rho\in(0,1)$. This feature adds a flexibility unprecedent in the aforementioned models. Besides the systematic component, the Unit-Weibull distribution is very flexible and can be regarded as an alternative to beta and Kumaraswamy distributions.

Inference in GARMA-like models is carried out typically via conditional \citep{Rocha2017, Bayerseas, Maior} or partial maximum likelihood estimation (PMLE) \citep{Pumi2017, BARC, Prass,helen}. We propose the use of the PMLE approach for parameter estimation because it is more general than the conditional approach and allows for the inclusion of random time-dependent covariates in the model. The PMLE approach also allows for the construction of confidence intervals, hypothesis testing and diagnostics, and forecasting future values. The paper is organized as follows. In the next section we introduce the proposed UWARMA model. In Section  \ref{infe} we discuss parameter estimation for UWARMA models presenting the score vector and  cumulative partial information matrix in closed form. Next we discuss inferential tools (Section \ref{inftools}), confidence interval, hypothesis testing, diagnostics (Section \ref{ht}) and forecasting (Section \ref{forecast}). A Monte Carlo simulation study is presented in Section \ref{MC} while in Section \ref{empi} we present an application of the UWARMA to a real dataset.

\section{The Unit-Weibull distribution and the UWARMA model}\label{model}

The Unit-Weibull distribution in its origins is obtained by the transformation $Y=e^{-X}$, for $X$ following a two-parameter Weibull distribution. The UW distribution, like the Kumaraswamy distribution, lacks a simple analytic expression for its expectation, which hinders its use in traditional mean-based regression models such as the beta regression of \cite{Ferrari2004}. However, the UW distribution has a simple closed formula for its quantile function. Taking advantage of this fact, \cite{UWreg} proposed a parameterization of the UW distribution in terms of its quantile function, allowing for the study of a quantile regression model similar to the Kumaraswamy regression of \cite{MK}. Under this parameterization, the UW distribution is absolutely continuous with respect to the Lebesgue measure in $(0,1)$ with density
\begin{equation}\label{dens}
f_\rho(y;\mu,\lambda):=\frac\lambda y\bigg[\frac{\log(\rho)}{\log(\mu)}\bigg]\bigg[\frac{\log(y)}{\log(\mu)}\bigg]^{\lambda-1}\rho^{\big[\frac{\log(y)}{\log(\mu)}\big]^\lambda}\!\!, \qquad y\in(0,1),
\end{equation}
where $\lambda>0$ is a shape parameter, $\mu\in(0,1)$ is the quantile parameter associated to the quantile $\rho\in(0,1)$ (assumed known) in the sense that, under \eqref{dens}, the $\rho$th quantile of the distribution is $\mu$. Thus, $\mu$ can be understood as a location parameter. If, for fixed $\rho\in(0,1)$, $Y$ is a random variable with density \eqref{dens}, we shall write $Y\sim \mathrm{UW}(\mu,\lambda;\rho)$. According to the values of its parameters, the UW distribution can present several shapes. For instance,  \eqref{dens} is unimodal when $\lambda>1$, bathtub shaped when $0<\lambda<1$, while for $\lambda=1$, \eqref{dens} presents J and inverse J shapes according to whether $\mu<\rho$ or $\mu>\rho$, respectively \citep{UWreg}.

Inspired by the Kumaraswamy regression model of \cite{MK} and following the ideas of the so-called generalized ARMA (GARMA) models of \cite{benjamin2003} and, in particular, the specification applied in the $\beta$ARMA models of \cite{Rocha2009}, \cite{Bayers} proposed the Kumaraswamy Autoregressive Moving Average (KARMA) models -  an observation driven model, where, conditionally to the process's history and a set of exogenous covariates, the random component follows a Kumaraswamy distribution reparameterized by its median, while the systematic component prescribes an ARMA-like structure for the conditional median in a generalized linear model (GLM) fashion, similar to the $\beta$ARMA's specification. Our goal is to apply the same approach to define a new dynamic time series model based on the UW distribution to model the conditional $\rho$th quantile through an ARMA-like structure.

Similar to the KARMA and $\beta$ARMA models, the proposed model is also an observation driven one. Let $\{Y_t\}_{t\in\Z}$ be a stochastic process taking values in $(0,1)$ and let $\{\bs x_t\}_{t\in\Z}$ be a set of $r$-dimensional exogenous covariates to be included in the model, which can be  either random or deterministic and time-dependent, or any combination of these. In observation driven models, the information available at time $t-1$ is denoted by $\F_{t-1}$ and usually is defined as a sigma-field of events. The flexibility regarding the nature of the covariates allowed in the present framework, however, creates a small nuisance when defining $\F_{t-1}$. On one hand, if all coordinates of $\bs x_t$ are random, then, at time $t-1$, the observer can only know its values up to time $t-1$. On the other hand, very commonly, at time $t-1$ the covariates are known up to time $t$ as, for instance, when they are deterministic functions of $t$ such as a polynomial trends, or when $\bs x_t$ itself is a lagged information related to another exogenous time-dependent process. In this case, it will be convenient to include this information in $\F_{t-1}$. When $\bs x_t$ is a combination of deterministic and random components, in which case, at time $t-1$, the deterministic components are known up to time $t$ while the random components are only known up to time $t-1$. To account for both cases, we write $\bs x_t=({\bs x^d_{t}}', {\bs x_{t-1}^s}')'$, where $\bs x^d_{t}$ denotes a $r_1$-dimensional non-random component, while $\bs x^s_{t}$ denotes a $r_2$-dimensional random component, at time $t$, with $0\leq r_1,r_2\leq r$, $r_1+r_2=r$. With this definition we  define $\F_{t-1}:=\sigma\{\bs x_t, Y_{t-1},\bs x_{t-1},Y_{t-2},\bs x_{t-2},\cdots\}$.

Let $\rho\in(0,1)$ be fixed. The model's random component is implicitly defined by prescribing that, conditionally on $\F_{t-1}$, $Y_t\sim\mathrm{UW}(\mu_t,\lambda;\rho)$ for some $\lambda>0$ and $\mu_t\in(0,1)$. In this condition, for all $t$, $\mu_t$ represents the conditional $\rho$th quantile of $Y_t$. Let $g:(0,1)\rightarrow\R$ be a twice differentiable strictly monotone link function. The systematic component of the proposed model is specified by
\begin{equation}\label{uwarma}
\eta_t:=g(\mu_t)= \alpha + \bs x_t^\prime \bs\beta + \sum_{i=1}^p \phi_i \big[ g(Y_{t-i})-\bs x_{t-i}^\prime \bs \beta\big] + \sum_{j=1}^q \theta_j r_{t-j},
\end{equation}
where $\eta_t$ is the linear predictor, $\alpha$ is an intercept, $\bs\beta=(\beta_1, \cdots,\beta_r)^\prime$ is the parameter vector related to the covariates, $\bs\phi=(\phi_1,\cdots,\phi_p)^\prime$ and $\bs\theta=(\theta_1,\cdots,\theta_q)^\prime$ are the AR and MA coefficients, respectively. The error term  in \eqref{uwarma} is defined in a recursive fashion just as in the $\beta$ARMA and KARMA models, namely, $r_t:=g(Y_t)-g(\mu_t)$. Observe that $\eta_t$ and $\mu_t$ are $\F_{t-1}$-measurable. Hereafter we shall assume that the following conditions, endemic to the discussion of classical ARMA models, hold:
\begin{enumerate}[(i)]
\item the AR and MA characteristic polynomials do not have common roots;
\item the AR characteristic polynomial do not have unit roots.
\end{enumerate}
Possible link functions to be applied in \eqref{uwarma} are the traditional logit, probit, loglog and cloglog, although parametric alternatives can also be considered, as in, for instance, \cite{Pumi2020}.

The proposed model, hereafter denoted UWARMA$(p,q)$, is defined by the specification $Y_t|\F_{t-1}\sim\mathrm{UW}(\mu_t,\lambda;\rho)$ and \eqref{uwarma}. Observe that \eqref{uwarma} has the same functional form as the $\beta$ARMA and KARMA models, although, unlike the KARMA model, proposed considering only a conditional median-based regression model, we define the model to allow any quantile to be explored. Explicit conditions for existence and stationarity of GARMA-like models can be obtained when the link function is the identity (i.e., when the support of the conditional distribution is $\R$), such as in the case of the SYMARMA/SYMARFIMA models \citep{Maior,helen}. For models designed for double-bounded time series, such as $\beta$ARMA and KARMA, these conditions are very hard to obtain and, to the best of our knowledge, apart from the specific work of \cite{BARC}, no general results are known.

\section{Parameter estimation}\label{infe}
In this section we consider parameter estimation in the context of UWARMA models. To accomplish that, we propose the use of partial maximum likelihood. This allow us to consider possibly time-dependent random and non-random covariates alike. Let $\rho\in(0,1)$ be fixed and let $Y_1,\cdots,Y_n$ be a sample from a UWARMA$(p,q)$ model with associated $r$-dimensional covariates $\bs x_1, \cdots, \bs x_n$. Let $\bs\gamma:=(\alpha,\bs\beta,\bs\phi,\bs\theta,\lambda)^\prime \in \Omega$, where $\Omega\subset\R^{r+p+q+1}\times(0,\infty)$ denotes the parameter space. Upon defining
$A_t:=\frac{\log(Y_t)}{\log(\mu_t)}$, let
\begin{equation*}
\ell_t(\bs\gamma):=\log(\lambda)-\log(Y_t)+\log\bigg(\frac{\log(\rho)}{\log(\mu_t)}\bigg)+(\lambda-1)\log(A_t)+\log(\rho)A_t^\lambda,
\end{equation*}
where $\mu_t$ is specified by \eqref{uwarma}. The partial log-likelihood is given by
\begin{equation}
\label{logL}
\ell(\bs \gamma)= \sum_{t=1}^{n} \ell_t(\bs\gamma),
\end{equation}
so that the partial maximum likelihood estimator (PMLE) of $\bs\theta$ is given by
\begin{equation*}\hat{\bs\gamma}=\underset{\bs\gamma\in\Omega}{\mathrm{argmax}}(\ell(\bs\gamma)).\end{equation*}
The partial score vector $\frac{\partial \ell(\bs\gamma)}{\partial \bs\gamma}$ can be used to obtain $\hat{\bs\gamma}$ by solving the system $\frac{\partial \ell(\bs\gamma)}{\partial \bs\gamma}=\bs 0$. In the next section,we derive closed formulas for the partial score vector, but the PMLE cannot be obtained analytically so that we have to resort to numerical optimization to accomplish that.

\subsection{Partial score vector}
By \eqref{logL}, in order to construct the partial score vector, it suffices to obtain the derivative of $\ell_t(\bs\gamma)$ with respect to $\bs\gamma$.  Let $\bs\nu:=(\alpha,\bs{\beta}^{\prime},\bs{\phi}^{\prime},\bs{\theta}^{\prime})^{\prime}$, so that $\bs\gamma = (\bs\nu',\lambda)'$. We start by computing
\begin{equation}\label{dldmu}
\frac{\partial \ell_t(\bs\gamma)}{\partial\lambda}= \frac1\lambda+\big( 1+ A_t^\lambda\log(\rho)\big)\log(A_t),
\quad\mbox{and}\quad
\frac{\partial \ell_t(\bs\gamma)}{\partial\mu_t} = -\frac{\lambda\big(1+\log(\rho)A_t^\lambda\big)}{\mu_t\log(\mu_t)}.
\end{equation}
Now, the chain rule gives
\begin{equation}\label{dldnu}
\frac{\partial \ell_t(\bs\gamma)}{\partial \nu_j} =\frac{\partial \ell_t(\bs\gamma)}{\partial\mu_t} \frac{\partial \mu_t}{\partial \eta_t}\frac{\partial \eta_t}{\partial \nu_j}=
 -\frac{\lambda\big(1+\log(\rho)A_t^\lambda\big)}{\mu_t\log(\mu_t)g'(\mu_t)}\bigg[\frac{\partial \eta_t}{\partial \nu_j}\bigg],
\end{equation}
where the last equality follows since $\frac{\partial\mu_t}{\partial \eta_t}=\frac1{g'(\mu_t)}$ and \eqref{dldmu}. Hence, we only need to obtain $\frac{\partial \eta_t}{\partial \nu_j}$. Since $\eta_t$ given in  \eqref{uwarma} is exactly the same specification for the KARMA model, the derivatives $\frac{\partial \eta_t}{\partial \nu_j}$ follow the same recursions as the ones presented in section 3 in \cite{Bayers}, namely
\begin{align}\label{equas}
\frac{\partial \eta_t}{\partial \alpha}&=1 - \sum_{j=1}^q \theta_j \frac{\partial \eta_{t-j}}{\partial \alpha},\qquad\qquad
\frac{\partial \eta_t}{\partial \beta_l}= x_{tl}-\sum_{i=1}^{p} \phi_ix_{(t-i)l}- \sum_{j=1}^q \theta_j \frac{\partial \eta_{t-j}}{\partial \beta_l},\nonumber\\
\frac{\partial \eta_t}{\partial \phi_k}&= g({y}_{t-1})-\bs x_{t-1}'\bs\beta - \sum_{j=1}^q \theta_j \frac{\partial \eta_{t-j}}{\partial \phi_k},\qquad\mbox{and}\qquad
\frac{\partial \eta_t}{\partial \theta_s}=r_{t-j} - \sum_{i=1}^q \theta_i \frac{\partial \eta_{t-i}}{\partial \theta_s},
\end{align}
for $l\in\{1,\cdots,r\}$, $k\in\{1,\cdots,p\}$ and $s\in\{1,\cdots,q\}$, where $x_{tl}$ denotes the $l$-th component of $\bs x_t$. Define
\begin{equation*}\bs h_1:=\bigg(\frac{\partial \ell_1(\bs\gamma)}{\partial \mu_1},\cdots, \frac{\partial \ell_n(\bs\gamma)}{\partial \mu_n}\bigg)^\prime \quad \mbox{and}\quad
\bs h_2:=\bigg(\frac{\partial \ell_1(\bs\gamma)}{\partial \lambda},\cdots, \frac{\partial \ell_n(\bs\gamma)}{\partial \lambda}\bigg)^\prime,
\end{equation*}
and let $D_{\bs\nu}$ be the $n\times(p+q+r+1)$ matrix with $(i,j)$th element given by
\begin{equation*}
[D_{\bs\nu}]_{i,j}:=\frac{\partial \ell_i(\bs\gamma)}{\partial \nu_j}.
\end{equation*}
The partial score vector $U(\bs\gamma)$ can be written as
\begin{equation*}
U(\bs\gamma)=\big(U_{\bs\rho}(\bs\gamma)',  U_\lambda(\bs\gamma)\big)', \qquad \mbox{with} \quad U_{\bs\nu}(\bs\gamma) := D_{\bs\nu}' \bs h_1 \quad\mbox{and}\quad  U_\lambda(\bs\gamma) := \bs1_n'\bs h_2,
\end{equation*}
where $\bs1_n:=(1,\cdots,1)'\in\R^n$.

\subsection{Partial information matrix}

In this section we derive the  information matrix per single information, denoted by $K(\boldsymbol{\bs\gamma})$. In the present context, direct knowledge of the unconditional distribution of the proposed model is also impossible to obtain. Hence, the traditional unconditional Fisher's information matrix is not obtainable. To obtain an analogous matrix in the context of dependent observations, we follow the ideas presented in \cite{Fokianos2004} and \cite{Kedem2002}. We start by defining the cumulative partial information matrix $K_n(\bs\gamma)$ by
\begin{equation*}K_n(\bs\gamma) := -\sum^n_{t=1}\E\left( \frac{\partial^2 \ell_t(\bs \gamma)}{\partial \bs\gamma\partial \bs \gamma^\prime } \Bigm| \mathscr{F}_{t-1}\right). \end{equation*}
We shall derive $K_n(\bs\gamma)$ in closed form. First notice that
\begin{align*}
\frac{\partial^2\ell_t(\bs\gamma)}{\partial \nu_i \partial \nu_j} &= \sum_{t=1}^{n}\frac{\partial}{\partial \mu_t}
\left( \frac{\partial \ell_t(\bs\gamma)}{\partial \mu_t}\frac{\partial \mu_t}{\partial \eta_t} \frac{\partial \eta_t}{\partial \nu_j}\right)
\frac{d \mu_t}{d \eta_t} \frac{\partial \eta_t}{\partial \nu_i} \\
&= \sum_{t=1}^{n} \left[ \frac{\partial^2 \ell_t(\bs\gamma)}{\partial \mu_t^2}\frac{\partial \mu_t}{\partial \eta_t} \frac{\partial \eta_t}{\partial \nu_j}
+ \frac{\partial \ell_t(\bs\gamma)}{\partial \mu_t}\frac{\partial}{\partial \mu_t}\left(\frac{\partial \mu_t}{\partial \eta_t} \frac{\partial \eta_t}{\partial \nu_j} \right) \right]
\frac{d \mu_t}{d \eta_t} \frac{\partial \eta_t}{\partial \nu_i}\,.
\end{align*}
From Lemma \ref{lema} in the Appendix, we conclude that $\E\Big(\frac{\partial \ell_t(\mu_t,\varphi)}{\partial \mu_t} \big| \F _{t-1}\Big)=0$ and, by the $\F_{t-1}$-mensurability of $\mu_t$ and $\eta_t$, it follows that
\begin{equation*}
\E\bigg(\frac{\partial^2\ell_t(\bs\gamma)}{\partial \nu_i \partial \nu_j}\Big|\F_{t-1}\bigg)=  \E\bigg(\frac{\partial^2 \ell_t(\bs\gamma)}{\partial \mu_t^2}\Big|\F_{t-1}\bigg)\bigg[\frac{\partial \mu_t}{\partial \eta_t}\bigg]^2 \frac{\partial \eta_t}{\partial \nu_i}\frac{\partial \eta_t}{\partial \nu_j},
\end{equation*}
where $\frac{\partial \eta_t}{\partial \nu_k}$ is given in \eqref{equas}. Elementary calculus yields
\begin{equation*}\label{d2ldmu2}
\frac{\partial^2 \ell_t(\bs\gamma)}{\partial \mu_t^2}=\frac{\lambda\big(\log(\mu_t)+1\big)\big(1+\log(\rho)A_t^\lambda\big)+\lambda^2\log(\rho)A_t^\lambda}{\mu_t^2\log(\mu_t)^2},
\end{equation*}
so that, by using Lemma \ref{lema}, we obtain
\begin{equation*}\E\bigg(\frac{\partial^2 \ell_t(\bs\gamma)}{\partial \mu_t^2}\Big|\F_{t-1}\bigg)=-\frac{\lambda^2}{\mu_t^2\log(\mu_t)^2}.\end{equation*}
Regarding the second derivatives with respect to $\lambda$, we have
\begin{equation*}\frac{\partial^2\ell_t(\bs\gamma)}{\partial\lambda^2}=\log(\rho)A_t^\lambda\log(A_t)^2-\frac1{\lambda^2},\end{equation*}
so that, by using \eqref{L2} and \eqref{L3} in Lemma \ref{lema}, we conclude that
\begin{align*}
%\begin{equation*}
\E\bigg(\frac{\partial^2 \ell_t(\bs\gamma)}{\partial \lambda^2}\Big|\F_{t-1}\bigg)&=
\frac{\pi^2(\kappa-2)\kappa-6\log\bigl(-\log(\rho)\big)\big[\log\bigl(-\log(\rho)\big)+2\kappa-2\big]}{6\lambda^2}-\frac1{\lambda^2},\\
&=-\frac1{\lambda^2}\bigg[1+\frac{\pi^2}{6}+(\kappa-2)\kappa+\log\bigl(-\log(\rho)\big)\big[\log\bigl(-\log(\rho)\big)\big]\bigg]\\
%&=\frac1{\lambda^2}\bigg[\frac{\pi^2(\kappa-2)\kappa}{6}-\log\bigl(-\log(\rho)\big)^2 -\big[2\kappa-2\big]\log\bigl(-\log(\rho)\big)-1\bigg]
%\end{equation*}
\end{align*}
where $\kappa=0.5772156649\dots$  is the Euler-Mascheroni constant \citep{GR2007}. Finally, since $\eta_t$ does not depend on $\lambda$, from \eqref{dldnu} we have
\begin{equation*}\frac{\partial^2\ell_t(\bs\gamma)}{\partial\lambda\partial\nu_j}=\frac\partial{\partial\lambda}\bigg(\frac{\partial \ell_t(\bs\gamma)}{\partial \nu_j}\bigg)=
-\frac{1+\log(\rho)A_t^\lambda+\lambda\log(\rho)A_t^\lambda\log(A_t)}{\mu_t\log(\mu_t)g'(\mu_t)}\bigg[\frac{\partial \eta_t}{\partial \nu_j}\bigg],\end{equation*}
and since $\eta_t$ is $\F_{t-1}$-measurable, the results in Lemma \ref{lema} yield
\begin{equation*}\E\bigg(\frac{\partial^2 \ell_t(\bs\gamma)}{\partial\lambda\partial\nu_j}\Big|\F_{t-1}\bigg)=\frac{1-\kappa-\log\big(-\log(\rho)\big)}{\mu_t\log(\mu_t)g'(\mu_t)} \bigg[\frac{\partial \eta_t}{\partial \nu_j}\bigg],\end{equation*}
with $\frac{\partial \eta_t}{\partial \nu_j}$ given in \eqref{equas}.
Let $T$ and $E_{\mu}$ be $n\times n$ diagonal matrices for which the $k$th diagonal elements are given by
\begin{equation*}[T]_{k,k}:=\frac1{g'(\mu_k)}\quad\mbox{and}\quad [E_\mu]_{k,k}:=-\E\bigg(\frac{\partial^2 \ell_k(\bs\gamma)}{\partial \mu_k^2}\Big|\F_{k-1}\bigg).\end{equation*}
Let $\bs {e}=(e_1,\cdots,e_n)'$ be the vector with $k$th coordinate given by
\begin{equation*}{e}_k:=-\E\bigg(\frac{\partial^2 \ell_k(\bs\gamma)}{\partial \mu_k\partial \lambda}\Big|\F_{k-1}\bigg)=\frac{\kappa+\log\big(-\log(\rho)\big)-1}{\mu_k\log(\mu_k)}. \end{equation*}
Upon writing
%$K_{\bs\nu,\bs\nu}:=D_{\bs\nu}'TE_{\mu}TD_{\bs\nu}$; $K_{\bs\nu,\lambda}=K_{\lambda,\bs\nu}'=D_{\bs\nu}'T\bs e$ and
\begin{equation*}
K_{\bs\nu,\bs\nu}:=D_{\bs\nu}'TE_{\mu}TD_{\bs\nu}, \quad  K_{\bs\nu,\lambda}=K_{\lambda,\bs\nu}'=D_{\bs\nu}'T\bs e\quad\mbox{and}\quad K_{\lambda,\lambda}:=\frac{1-2\big[\kappa+\log\big(-\log(\rho)\big)\big]}\lambda+\frac1{\lambda^2}
\end{equation*}
we have
\begin{equation*}
K_n(\bs\gamma)=\left(
\begin{array}{cc}
  K_{\bs\nu,\bs\nu} & K_{\bs\nu,\lambda} \\
  K_{\lambda,\bs\nu} & K_{\lambda,\lambda}
\end{array}
\right).
\end{equation*}
Under mild conditions, it can be shown that $\frac1nK_n(\bs\gamma){\ \longrightarrow \ } K(\bs\gamma)$ in probability, where $K(\bs\gamma)$ is a positive definite matrix, called information matrix per single information \citep{Fokianos2004}. This is the analogous of the Fisher's information matrix per single information in the unconditional case.

\section{Inferential tools}\label{inftools}

Under very mild conditions, the asymptotic theory of the PMLE for general GARMA-like models is derived in \cite{Fokianos1998, Fokianos2004}. However, the aforementioned works only consider the case where the random component is a member of the canonical exponential family. As such, the estimation of the distributional parameter (associated to the conditional density)  and the parameters associated to the systematic components are independent, so that the authors derive the asymptotic theory for the PMLE assuming that the distribution parameter is known.

For distributions which are not members of the canonical exponential family, usually the estimation of the distribution parameter is not independent of the other parameters. This is the case, for instance, for the $\beta$ARMA, KARMA and the proposed UWARMA. In these situations, the arguments presented in \cite{Fokianos1998, Fokianos2004} seem applicable in a case-by-case fashion, under the same conditions, with the necessary adaptations and by considering the distribution parameter estimated along with the other parameters. Details will be left for future research.

An important facet of the theory presented in \cite{Fokianos1998, Fokianos2004} is that the assumptions related to the model's systematic component are high level ones, which are very hard to translate into assumptions related to the dynamic of $\beta$ARMA, KARMA and UWARMA alike, and generally difficult to verify. Nevertheless, the assumptions are very mild, mainly intending to assure that the model is well-defined and that $\frac{K_n(\bs\gamma)}n\rightarrow K(\bs \gamma)$ in probability, with $K(\bs\gamma)$ positive definite.

Let $Y_1, \cdots, Y_n$  be a sample from a UWARMA$(p,q)$ model with with true parameter $\bs\gamma_0$ satisfying conditions (i) e (ii) in Section \ref{model}, and let $\bs x_1, \cdots, \bs x_n$ be a set of $r$-dimensional, possibly random and time-dependent, exogenous covariates to be included in the model. Let $\widehat{\bs\gamma}$ denote a solution of $U(\bs\gamma)=\bs0$. Under mild conditions \citep[closely related to the ones presented in][]{Fokianos2004}, the PMLE will be consistent,
\[\widehat{\bs\gamma}_n\overset{P}{\underset{n\rightarrow\infty}\longrightarrow} \bs\gamma_0,\]
and asymptotically normal
\begin{equation}\label{an}
\sqrt{n}(\widehat{\bs\gamma}_n-\bs\gamma_0)\overset{d}{\underset{n\rightarrow\infty}\longrightarrow} N_{p+q+r+1}\big(\bs0,K(\bs\gamma_0)^{-1}\big).
\end{equation}

\section{Hypothesis tests, confidence intervals and diagnostics}\label{ht}
Construction of asymptotic tests and confidence intervals for the proposed UWARMA model can be carried on in the same lines as those presented for the KARMA and $\beta$ARMA models. Considering the framework of the previous section, let $\hat\gamma_j$ denote the $j$th component of the PMLE $\hat{\bs\gamma}$ based on a sample of size $n$ from a UWARMA$(p,q)$ model, respectively, and assume that \eqref{an} holds. In this case, we can easily construct asymptotic confidence intervals and perform asymptotic hypothesis testing. Observe that from \eqref{an}, for large $n$,
\begin{equation}\label{approx}
\big[K_n(\hat{\bs\gamma})^{jj}\big]^{-\frac12}(\hat\gamma_j-\gamma_{0j})\approx N(0,1),
\end{equation}
where $\gamma_{0j}$ denotes the $j$th component of the true parameter $\bs\gamma_0$ and $K_n(\hat{\bs\gamma})^{jj}$ denotes the $j$th diagonal element of $K_n(\hat{\bs\gamma})^{-1}$. From \eqref{approx}, for $0<\delta<1/2$, a $100(1-\delta)\%$ asymptotic confidence interval for $\gamma_j$ is given by
\begin{equation*}
\widehat{\gamma}_j \pm z_{1-\delta/2} \big(K_n(\widehat{\boldsymbol{\gamma}})^{jj}\big)^{1/2},
\end{equation*}
where $z_{1-\delta/2}$ denotes the $1-\delta/2$ quantile of the standard normal distribution. Asymptotic test statistics for commonly applied tests, such as Wald's $z$, Rao's score, likelihood ratio, etc. can be derived from \eqref{an} as well. For instance, testing hypothesis of the form $H_0:\gamma_j=\gamma_j^\ast$ against $H_1:\gamma_j\neq\gamma_j^\ast$ for some $\gamma_j^\ast$ given,  can be done applying \eqref{approx} and the traditional Wald's $z$ statistics, namely
\begin{equation*}
z=\frac{\hat{\gamma}_j-{\gamma}_j^\ast}{\big(K_n(\widehat{\boldsymbol{\gamma}})^{jj}\big)^{1/2}},
\end{equation*}
which, under $H_0$ and for large enough $n$ is approximately distributed as $N(0,1)$.  Other tests follow the same principle and, for large samples, their null distribution will be the same as their counterparts based on independent samples.

Information criteria can be useful for automatic model comparison. In the context of UWARMA,  Akaike's AIC, Schwartz's BIC and Hannan-Quinn's HQC are obtained by considering the maximized partial log-likelihood $\ell(\hat{\bs\gamma})$ and analyzed as per usual. Portmanteau tests can also be defined in the usual fashion \citep{port}. Finally, we mention that residual analysis in the lines of \cite{Rocha2009} for the $\beta$ARMA and \cite{Bayers} for the KARMA, can also be performed in the context of UWARMA models. Observe, however, that we make no assumption regarding the UWARMA's residual, so that residuals are but tools to further explore the model, and not requirements for its validity.

\section{Forecast}\label{forecast}
In-sample and out-of-sample forecast in the context of UWARMA models follow the same recipe as in the case of $\beta$ARMA and KARMA models. Let $\hat{\bs\gamma}$ be the PMLE based on a sample $y_1,\cdots,y_n$ from a UWARMA$(p,q)$ with associated covariates $\bs x_1,\cdots,\bs x_n$. In-sample forecasts are obtained by recursively reconstructing a sequence $\hat{\mu}_1, \cdots,\hat{\mu}_n$ through \eqref{uwarma}. Observe that in this case we are using information provided by the $\rho$th quantile to forecast the time series. From $\hat{\mu}_t$ we can also reconstruct the error term by setting $\hat{r}_t:=\big(g(y_t)-g(\hat\mu_t)\big)I(0\leq t\leq n)$. More precisely, we write
\[\widehat\mu_t=g^{-1} \bigg( \widehat{\alpha}+ \bs x_t'\widehat{\bs\beta} + \sum_{i=1}^{p}\widehat{\phi}_i \big( g(y_{t-i})- \bs x_{t-i}'\widehat{\bs\beta}\big) +\sum_{j=1}^{q}\widehat{\theta}_j \widehat{r}_{t-j} \bigg)I(0\leq t\leq n).\]
Similarly, $h$-step ahead forecasts, say $\hat y_{n+1}, \cdots, \hat y_{n+h}$ can be obtained in a similar fashion. We define
\begin{equation}\label{for}
\hat{y}_{n+h}:=g^{-1} \bigg( \widehat{\alpha}+ \bs x_{n+h}'\widehat{\bs\beta} + \sum_{i=1}^{p}\widehat{\phi}_i \big(g([y_{n+h-i}]^\ast)- \bs x_{n+h-i}'\widehat{\bs\beta}\big) +\sum_{j=1}^{q}\widehat{\theta}_j \widehat{r}_{n+h-j} \bigg),
\end{equation}
where $[y_t]^\ast:= y_tI(1\leq t \leq n)+\hat{y}_tI(t\geq n+1)$. Observe that in the presence of covariates, \eqref{for} tacitly requires $h$-step ahead values for the covariates to be provided. This is very easy when the covariates are deterministic, such as in the case of a polynomial trend or when a sine/cosine function is applied to model a deterministic seasonality. If the covariates are themselves time series, one may need to model and forecast the covariates before proceeding with the forecasting exercise.

\section{Monte Carlo Simulation}\label{MC}

In this section we present a Monte Carlo simulation study to evaluate the finite sample and forecast performance of the PMLE in the context of the proposed UWARMA models. To generate samples from a UWARMA$(p,q)$ process with parameter $\bs\gamma$ for $\rho\in(0,1)$ fixed and $\mu_t$ given by \eqref{uwarma}, by the $\F_{t-1}$-measurability of $\mu_t$, we only need to generate independent random variates $y_t$ with distribution UW$(\mu_t,\lambda;\rho)$. Observe that the dependence in the time series is induced by the dependence in the sequence $\mu_t$ which is the $\rho$th quantile of $Y_t$. For a given $\mu\in(0,1)$, generating a random variate $y$ from  $Y\sim\mathrm{UW}(\mu,\lambda;\rho)$ is straightforward: simply generate $u\sim U(0,1)$ and take
\[y:=\mu e^{\big[\frac{\log(1-u)}{\log(\rho)}\big]^{\frac1\lambda}}.\]%\mu_t\exp\bigg\{\frac{\log(1-u)}{\log(\rho)}\bigg\}^{\frac1\lambda}=
To generate a sample from a UWARMA$(p,q)$ model, the idea is iteratively construct $\mu_t$ and sample $y_t$. For $t\leq 0$, we initialize $y_t$,  $\mu_t$, $r_t$ with 0 and $\bs x_t=\bs 0$. Then we recursively construct $\mu_t$ for $t\geq1$. For instance, for $t=1$ we have $\mu_1=g^{-1}\big(\alpha+\bs x_1'\bs\beta\big)$, $r_1=g(y_1)-g(\mu_1)$ and finally $y_1\sim \mathrm{UW}(\mu_1,\lambda;\rho)$; for $t=2$ we have
$\mu_2=g^{-1}\big(\alpha+\bs x_2'\bs\beta+\phi_1\big[g(y_{1})-\bs x_{1}^\prime\bs \beta\big]+\theta_1 r_1\big)$, $r_2=g(y_2)-g(\mu_2)$ and $y_2\sim \mathrm{UW}(\mu_2,\lambda;\rho)$ and so on.

\begin{figure}[h!]
\centering
\includegraphics[width=\textwidth]{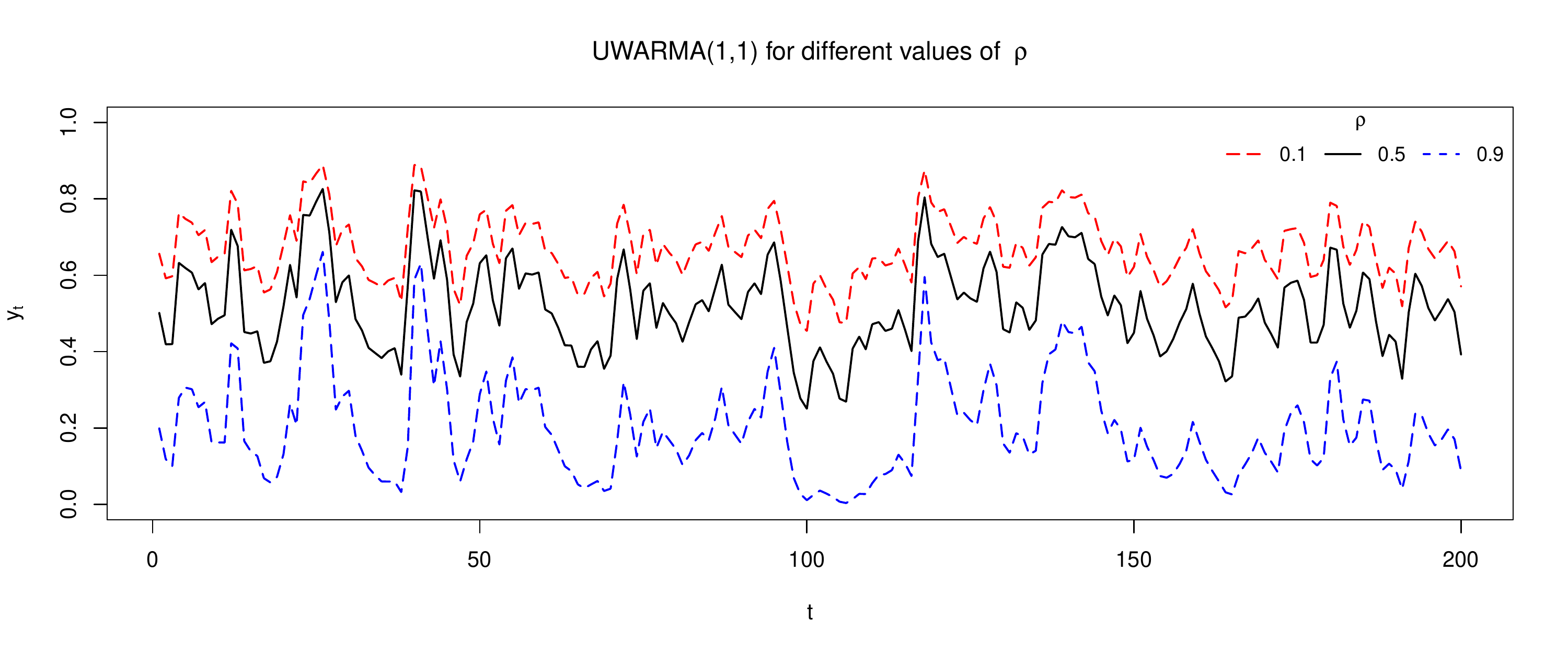}
\caption{Simulated UWARMA$(1,1)$ with $\lambda = 6$, $\phi=0.4$, $\theta=0.6$ and $\rho\in\{0.1,0.5,0.9\}$, produced using the same random seed.} \label{exe}
\end{figure}

Figure \ref{exe} presents the simulated sample path of a UWARMA$(1,1)$ with $\lambda = 6$, $\phi=0.4$, $\theta=0.6$, $\alpha=0$, no other covariates and $\rho\in\{0.1,0.5,0.9\}$. The sample paths were generated using the same random seed.
Figure \ref{exe} shows a facet of UWARMA's with same structure but different $\rho$'s, which may seem counterintuitive at first glance: for small values of $\rho$ a sample path of an UWARMA tends to be concentrated on the upper side of $(0,1)$ while for higher values, the opposite happens. To understand why this is the case, suppose that we have an UWARMA$(p,q)$ process $Y_t$ with $\alpha=0$ and no other covariates, and suppose $\rho=0.1$. Suppose that, for some $t_0$, we have $\mu_{t_0}=1/2$ which means that $\eta_{t_0}=0$. Since $\mu_t$ is the conditional $\rho$th quantile of $Y_t$, the value of $Y_{t_0}\sim UW(\mu_{t_0},\lambda;0.1)$ will be higher than 1/2 with high probability. Given the structure of $\eta_t$, $\mu_t$ tends to oscillate around 1/2 and, hence, a sample path of $Y_t$ will most often assume values greater than 1/2. On the other hand, suppose that $Z_t$ is an UWARMA$(p,q)$ with similar structure as $Y_t$ but with $\rho=0.9$. Suppose again that $\mu_{t_0}=1/2$. In this case the value of $Z_{t_0}\sim UW(\mu_{t_0},\lambda;0.9)$ will very likely to be smaller than 1/2 and, for the same reason stated before, $Z_t$ will likely to assume values smaller than 1/2 most of the time. Therefore for UWARMA's with the same structure, a sample path for which $\rho=\rho_0$ is likely to be above the one for which $\rho>\rho_0$.

Since $\eta_t$ has an ARMA-like structure, it is likely to stay in a neighborhood of 0, so that $\eta_{t_0+1}$ will likely be close to 0 and $\mu_{t_0+1}$ is likely to be close to 1/2. To fix the ideas, that this is the case.
Since $\mu_t$ is the conditional $\rho$th quantile of $Y_t$, if $\rho = 0.1$, say, then we expect that with high probability, the value of $Y_{t_0}\sim UW(\mu_{t_0},\lambda;0.1)$ will likely to be above 1/2, since $\rho=0.1$ implies that 90\% of the values observed for $Y_{t_0}$

In this type of model, it is advisable to perform a small burn-in to counter any effect of initialization.

Another noteworthy facet of generating GARMA-like models as well such as KARMA, $\beta$ARMA and UWARMA is the numerically instability some combination of parameters may cause. The literature on this matter is remarkably vague. One of the few works that acknowledges numerically instabilities in simulating $\beta$ARMA models is \cite{Casarin2012}. Numerical instability for UWARMA (and GARMA models taking values in $(0,1)$ in general) occurs mainly in two scenarios: when $|\eta_t|$ is large so that $\mu_t = g^{-1}(\eta_t)$ is so close to the boundaries of $(0,1)$ beyond the computer's capabilities to distinguish it from 0 or 1; and when small neighborhoods of 0 and 1 present high probability of occurrence, which causes the sample path to travel too close to the extremes often, eventually becoming numerically indistinguishable from 0 or 1. For the UW distribution, this happens when $\lambda<1$ in which case the density is bathtube-shaped presenting asymptotes in the extremes \citep{UWreg}. This makes simulating samples from UWARMA models with $\lambda<1$ extremely difficult.

Regarding inference, numerical optimization of the partial log-likelihood requires initialization. In the simulation, initial values for $(\alpha,\bs\beta',\bs\phi')$ are given by the ordinary least square estimate of the regression
\[g(y_t)= \alpha + \bs x_t^\prime \bs\beta + \sum_{i=1}^p \phi_i g(y_{t-i})+\varepsilon_t,\]
where $\varepsilon_t$ is a generic error term. We initialize $\lambda=10$ and $\bs\theta=\bs0$.

\subsection*{DGP}

We simulate 1,000 replicas of the UWARMA(1,1) model with $\rho\in\{0.25,0.5,0.75\}$, $\lambda\in\{5,10,20\}$, $\alpha=0$ and $(\phi,\theta)\in\{(0.6,0.4),(0.4,0.6),(-0.4,0.6),(0.4,-0.6)\}$ and sample sizes $n\in\{250,500,1000\}$. All simulations were performed in R \citep{R} version 4.1.3. To simulate and fit the proposed UWARMA, we use package \texttt{BTSR}. A burn-in of size 1,000 was applied in all simulations.

\subsection*{Results}
Simulation results are presented in Table \ref{rho25}, with the exception of $\alpha$ whose estimated values were omitted for presentation simplicity. From the table we observe that the PMLE performs very well with small biases even for sample size 250. Overall the estimation variability is small and decreases with $n$. As expected, $\lambda$ present the highest variability. The worst result happens for the combination $\lambda = 20$, $\phi=-0.4$, $\theta = 0.6$ and $\rho=0.75$, for which the estimates for $\phi$ and $\theta$ present high bias for all sample sizes due to a flat likelihood surface for this particular combination.
\begin{table}
\caption{Simulation Results. Presented are the estimated values and standard deviation (in parenthesis). Estimated values of $\alpha$ were omitted to save space. }\label{rho25}
\centering
\scriptsize
\renewcommand{\arraystretch}{1.1}
\setlength{\tabcolsep}{4pt}
\vspace{.3cm}
\begin{tabular}{c|c|ccc|ccc|ccc|ccc}
\hline
\multirow{2}{*}{$\rho$}&$\multirow{2}{*}{$\lambda$}$&\multicolumn{3}{c|}{$(\phi,\theta)=(0.6,0.4)$}&\multicolumn{3}{c|}{$(\phi,\theta)=(0.4,0.6)$}&
\multicolumn{3}{c|}{$(\phi,\theta)=(-0.4,0.6)$}&\multicolumn{3}{c}{$(\phi,\theta)=(0.4,-0.6)$}\\
\cline{3-14}
&&$\hat\phi$&$\hat\theta$&$\hat\lambda$&$\hat\phi$&$\hat\theta$&$\hat\lambda$&$\hat\phi$&$\hat\theta$&$\hat\lambda$&$\hat\phi$&$\hat\theta$&
$\hat\lambda$\\
\hline
\multicolumn{2}{c}{}&\multicolumn{12}{c}{$n=250$}\\
\hline
\multirow{6}{*}{\begin{sideways}0.25\end{sideways}}
&\multirow{2}{*}{5} &0.591  &    0.389  &    5.010  &    0.388  &    0.591  &    5.003  &    -0.377  &    0.578  &    5.061  &    0.363  &    -0.567  &    5.049\\
                    & & (0.054)& (0.065)& (0.262)& (0.066)& (0.058)& (0.265)& (0.232)& (0.220)& (0.262)& (0.252)& (0.243)& (0.262) \\
&\multirow{2}{*}{10}&0.591  &    0.395  &    10.016  &    0.391  &    0.588  &    9.987  &    -0.366  &    0.571  &    10.121  &    0.352  &    -0.555  &    10.113\\
                    & & (0.053)& (0.062)& (0.517)& (0.064)& (0.054)& (0.524)& (0.218)& (0.205)& (0.514)& (0.250)& (0.247)& (0.510) \\
&\multirow{2}{*}{20}&0.586  &    0.396  &    19.938  &    0.383  &    0.596  &    19.938  &    -0.370  &    0.573  &    20.263  &    0.221  &    -0.432  &    20.067\\
                    & & (0.058)& (0.070)& (1.361)& (0.073)& (0.061)& (1.219)& (0.193)& (0.177)& (1.013)& (0.266)& (0.273)& (1.163) \\
\hline
\multirow{6}{*}{\begin{sideways}0.5\end{sideways}}
&\multirow{2}{*}{5} &0.591  &    0.403  &    5.032  &    0.394  &    0.598  &    5.020  &    -0.371  &    0.574  &    5.062  &    0.383  &    -0.589  &    5.060\\
                    & & (0.053)& (0.057)& (0.250)& (0.065)& (0.058)& (0.253)& (0.219)& (0.206)& (0.252)& (0.225)& (0.214)& (0.256) \\
&\multirow{2}{*}{10}&0.597  &    0.397  &    10.044  &    0.398  &    0.596  &    10.000  &    -0.356  &    0.561  &    10.115  &    0.375  &    -0.580  &    10.140\\
                    & & (0.052)& (0.058)& (0.508)& (0.059)& (0.060)& (0.658)& (0.222)& (0.209)& (0.521)& (0.226)& (0.215)& (0.525) \\
&\multirow{2}{*}{20}&0.591  &    0.401  &    20.126  &    0.390  &    0.601  &    20.073  &    -0.321  &    0.527  &    20.064  &    0.283  &    -0.489  &    20.129\\
                    & & (0.051)& (0.058)& (1.037)& (0.066)& (0.054)& (1.159)& (0.200)& (0.187)& (1.464)& (0.233)& (0.233)& (1.109) \\
\hline
\multirow{6}{*}{\begin{sideways}0.75\end{sideways}}
&\multirow{2}{*}{5} & 0.547  &    0.418  &    4.845  &    0.372  &    0.604  &    4.986  &    -0.367  &    0.572  &    5.05  &    0.369  &    -0.573  &    5.054\\
                    && (0.093)& (0.061)& (0.455)& (0.066)& (0.054)& (0.297)& (0.200)& (0.187)& (0.258)& (0.224)& (0.213)& (0.254) \\
&\multirow{2}{*}{10}& 0.551  &    0.416  &    9.765  &    0.372  &    0.604  &    9.970  &    -0.354  &    0.559  &    10.117  &    0.363  &    -0.572  &    10.139\\
                    && (0.087)& (0.061)& (0.755)& (0.067)& (0.052)& (0.665)& (0.215)& (0.204)& (0.531)& (0.224)& (0.212)& (0.515) \\
&\multirow{2}{*}{20}& 0.559  &    0.408  &    19.489  &    0.381  &    0.590  &    19.482  &    -0.157  &    0.354  &    19.977  &    0.386  &    -0.592  &    20.218\\
                    && (0.081)& (0.074)& (1.886)& (0.072)& (0.073)& (2.047)& (0.288)& (0.303)& (1.362)& (0.189)& (0.172)& (1.079) \\
\hline
\multicolumn{2}{c}{}&\multicolumn{12}{c}{$n=500$}\\
\hline
\multirow{6}{*}{\begin{sideways}0.25\end{sideways}}
&\multirow{2}{*}{5} & 0.593  &    0.395  &    5.009  &    0.395  &    0.593  &    5.001  &    -0.386  &    0.589  &    5.031  &    0.385  &    -0.586  &    5.023\\
                    && (0.039)& (0.044)& (0.181)& (0.045)& (0.042)& (0.188)& (0.139)& (0.126)& (0.175)& (0.146)& (0.131)& (0.195) \\
&\multirow{2}{*}{10}& 0.594  &    0.393  &    9.988  &    0.394  &    0.595  &    9.982  &    -0.383  &    0.584  &    10.053  &    0.375  &    -0.578  &    10.050\\
                    && (0.037)& (0.044)& (0.347)& (0.045)& (0.041)& (0.423)& (0.131)& (0.119)& (0.368)& (0.146)& (0.129)& (0.374) \\
&\multirow{2}{*}{20}& 0.592  &    0.399  &    19.884  &    0.390  &    0.595  &    19.911  &    -0.385  &    0.585  &    20.115  &    0.275  &    -0.481  &    19.922\\
                    && (0.041)& (0.053)& (1.221)& (0.075)& (0.045)& (1.085)& (0.123)& (0.111)& (0.730)& (0.202)& (0.196)& (0.895) \\
\hline
\multirow{6}{*}{\begin{sideways}0.5\end{sideways}}
&\multirow{2}{*}{5} &  0.598  &    0.398  &    5.006  &    0.395  &    0.597  &    5.006  &    -0.381  &    0.583  &    5.025  &    0.392  &    -0.595  &    5.022\\
                    &&  (0.036)& (0.040)& (0.181)& (0.043)& (0.041)& (0.192)& (0.144)& (0.131)& (0.179)& (0.147)& (0.134)& (0.172) \\
&\multirow{2}{*}{10}&  0.598  &    0.398  &    10.024  &    0.401  &    0.596  &    9.990  &    -0.384  &    0.585  &    10.078  &    0.374  &    -0.577  &    10.051\\
                    &&  (0.036)& (0.041)& (0.361)& (0.043)& (0.045)& (0.512)& (0.134)& (0.123)& (0.367)& (0.145)& (0.132)& (0.346) \\
&\multirow{2}{*}{20}&  0.594  &    0.402  &    20.046  &    0.398  &    0.597  &    20.011  &    -0.349  &    0.550  &    19.897  &    0.325  &    -0.525  &    20.124\\
                    &&  (0.043)& (0.045)& (0.779)& (0.051)& (0.039)& (0.831)& (0.152)& (0.140)& (1.223)& (0.168)& (0.164)& (0.857) \\
\hline
\multirow{6}{*}{\begin{sideways}0.75\end{sideways}}
&\multirow{2}{*}{5} &0.566  &    0.412  &    4.867  &    0.388  &    0.598  &    4.971  &    -0.397  &    0.597  &    5.034  &    0.395  &    -0.594  &    5.026\\
                    & & (0.062)& (0.042)& (0.346)& (0.046)& (0.043)& (0.218)& (0.132)& (0.117)& (0.174)& (0.141)& (0.129)& (0.179) \\
&\multirow{2}{*}{10}&0.577  &    0.408  &    9.820  &    0.385  &    0.602  &    9.956  &    -0.363  &    0.564  &    10.062  &    0.392  &    -0.594  &    10.045\\
                    & & (0.051)& (0.040)& (0.533)& (0.044)& (0.036)& (0.460)& (0.151)& (0.142)& (0.385)& (0.141)& (0.127)& (0.349) \\
&\multirow{2}{*}{20}&0.578  &    0.407  &    19.527  &    0.393  &    0.594  &    19.586  &    -0.152  &    0.350  &    19.800  &    0.386  &    -0.590  &    20.061\\
                    & & (0.049)& (0.042)& (1.416)& (0.049)& (0.049)& (1.722)& (0.253)& (0.265)& (1.125)& (0.134)& (0.120)& (0.849) \\
\hline

\multicolumn{2}{c}{}&\multicolumn{12}{c}{$n=1{,}000$}\\
\hline
\multirow{6}{*}{\begin{sideways}0.25\end{sideways}}
&\multirow{2}{*}{5} & 0.599  &    0.395  &    5.007  &    0.398  &    0.594  &    5.000  &    -0.394  &    0.597  &    5.010  &    0.390  &    -0.587  &    5.012\\
                    && (0.027)& (0.030)& (0.122)& (0.032)& (0.036)& (0.134)& (0.093)& (0.082)& (0.121)& (0.100)& (0.090)& (0.126) \\
&\multirow{2}{*}{10}& 0.597  &    0.396  &    10.000  &    0.398  &    0.596  &    9.978  &    -0.385  &    0.586  &    10.015  &    0.380  &    -0.581  &    10.018\\
                    && (0.026)& (0.032)& (0.254)& (0.033)& (0.030)& (0.340)& (0.095)& (0.085)& (0.256)& (0.105)& (0.097)& (0.263) \\
&\multirow{2}{*}{20}& 0.595  &    0.402  &    19.804  &    0.389  &    0.598  &    19.881  &    -0.389  &    0.589  &    20.063  &    0.301  &    -0.503  &    19.913\\
                    && (0.032)& (0.049)& (1.200)& (0.073)& (0.032)& (0.826)& (0.089)& (0.079)& (0.518)& (0.185)& (0.182)& (0.673) \\
\hline
\multirow{6}{*}{\begin{sideways}0.5\end{sideways}}
&\multirow{2}{*}{5} & 0.600  &    0.400  &    5.007  &    0.403  &    0.592  &    5.005  &    -0.394  &    0.594  &    5.011  &    0.394  &    -0.594  &    5.010\\
                    && (0.026)& (0.029)& (0.127)& (0.039)& (0.059)& (0.128)& (0.093)& (0.085)& (0.129)& (0.097)& (0.088)& (0.120) \\
&\multirow{2}{*}{10}& 0.599  &    0.399  &    10.016  &    0.403  &    0.590  &    9.852  &    -0.389  &    0.589  &    10.024  &    0.378  &    -0.580  &    10.040\\
                    && (0.025)& (0.028)& (0.252)& (0.040)& (0.056)& (0.840)& (0.094)& (0.087)& (0.254)& (0.116)& (0.104)& (0.280) \\
&\multirow{2}{*}{20}& 0.598  &    0.400  &    20.024  &    0.398  &    0.599  &    19.997  &    -0.361  &    0.563  &    19.932  &    0.330  &    -0.531  &    20.066\\
                    && (0.029)& (0.029)& (0.531)& (0.030)& (0.026)& (0.506)& (0.118)& (0.109)& (0.920)& (0.147)& (0.146)& (0.659) \\
\hline
\multirow{6}{*}{\begin{sideways}0.75\end{sideways}}
&\multirow{2}{*}{5} &  0.583  &    0.407  &    4.918  &    0.405  &    0.575  &    4.976  &    -0.393  &    0.595  &    5.011  &    0.389  &    -0.588  &    5.011\\
                    &&  (0.040)& (0.026)& (0.249)& (0.052)& (0.088)& (0.151)& (0.095)& (0.084)& (0.123)& (0.089)& (0.084)& (0.121) \\
&\multirow{2}{*}{10}&  0.586  &    0.405  &    9.894  &    0.396  &    0.597  &    9.923  &    -0.363  &    0.564  &    10.048  &    0.388  &    -0.590  &    10.043\\
                    &&  (0.033)& (0.030)& (0.402)& (0.034)& (0.030)& (0.539)& (0.116)& (0.114)& (0.292)& (0.095)& (0.089)& (0.271) \\
&\multirow{2}{*}{20}&  0.586  &    0.405  &    19.683  &    0.396  &    0.595  &    19.632  &    -0.168  &    0.362  &    19.643  &    0.392  &    -0.594  &    20.070\\
                    &&  (0.035)& (0.029)& (1.202)& (0.037)& (0.039)& (1.704)& (0.255)& (0.268)& (1.309)& (0.090)& (0.078)& (0.525) \\
\hline
\end{tabular}
\end{table}

\subsection*{Joint behavior}

\begin{figure}[h!]
\centering
\mbox{
\includegraphics[width=.5\textwidth]{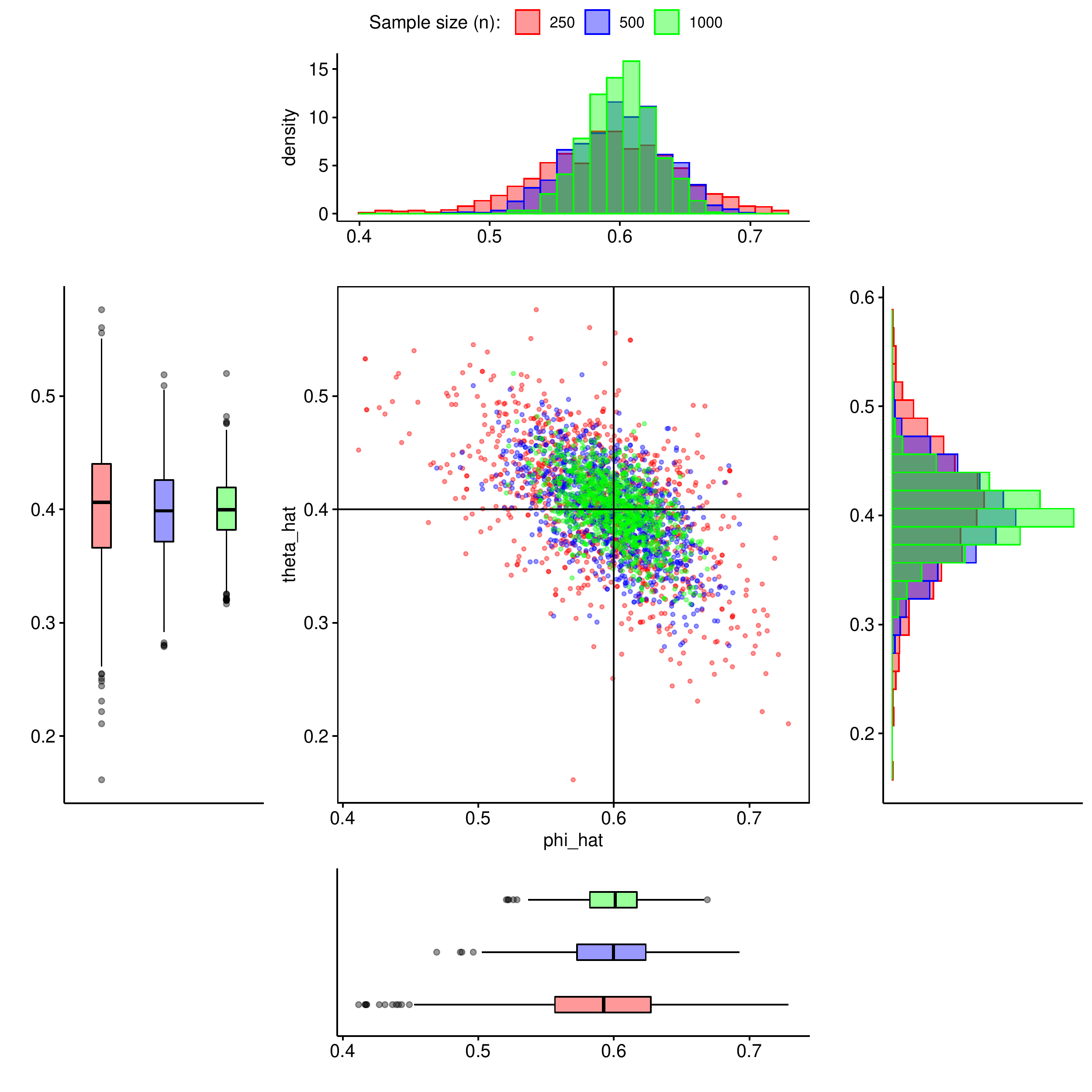}
\includegraphics[width=.5\textwidth]{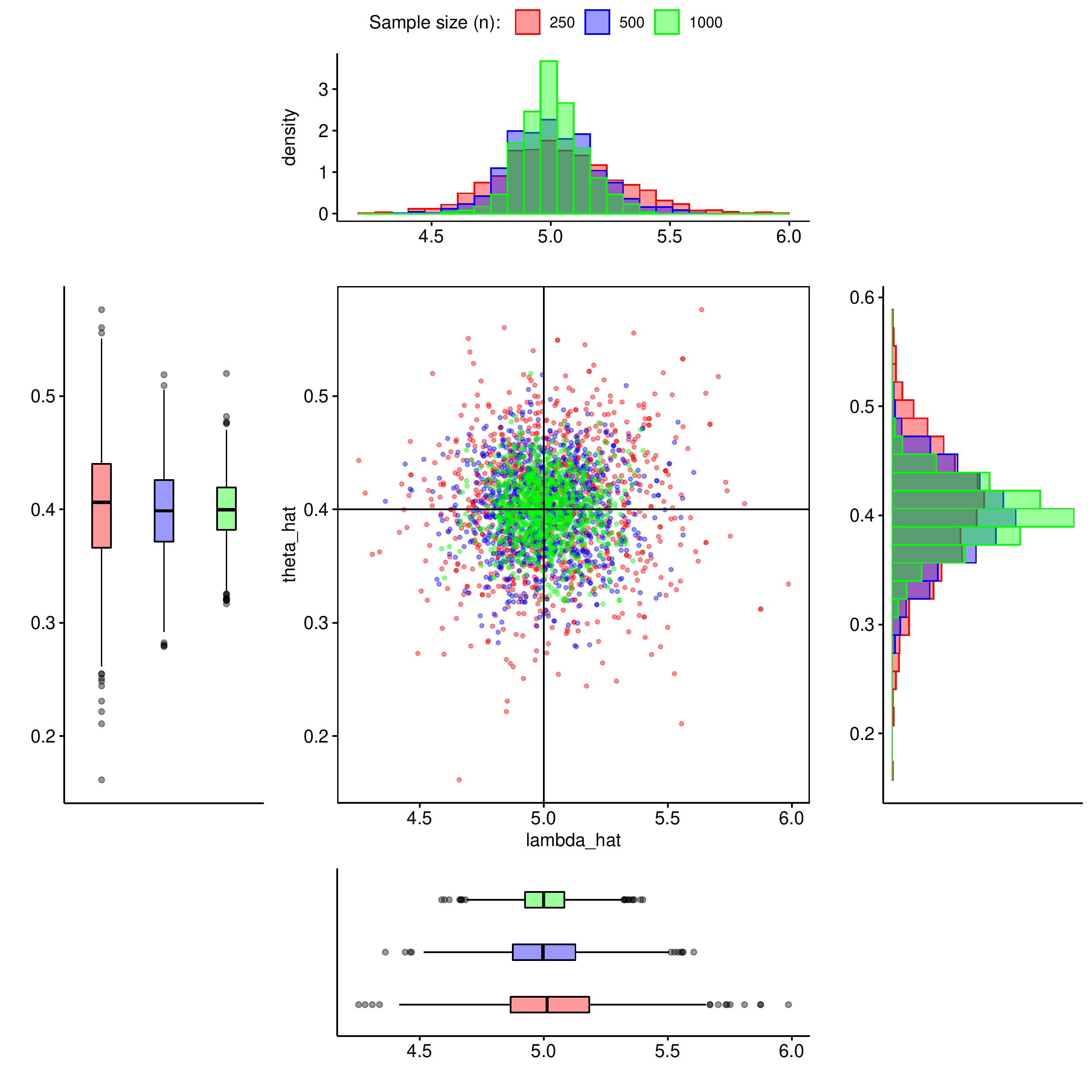}
}
\mbox{\includegraphics[width=.5\textwidth]{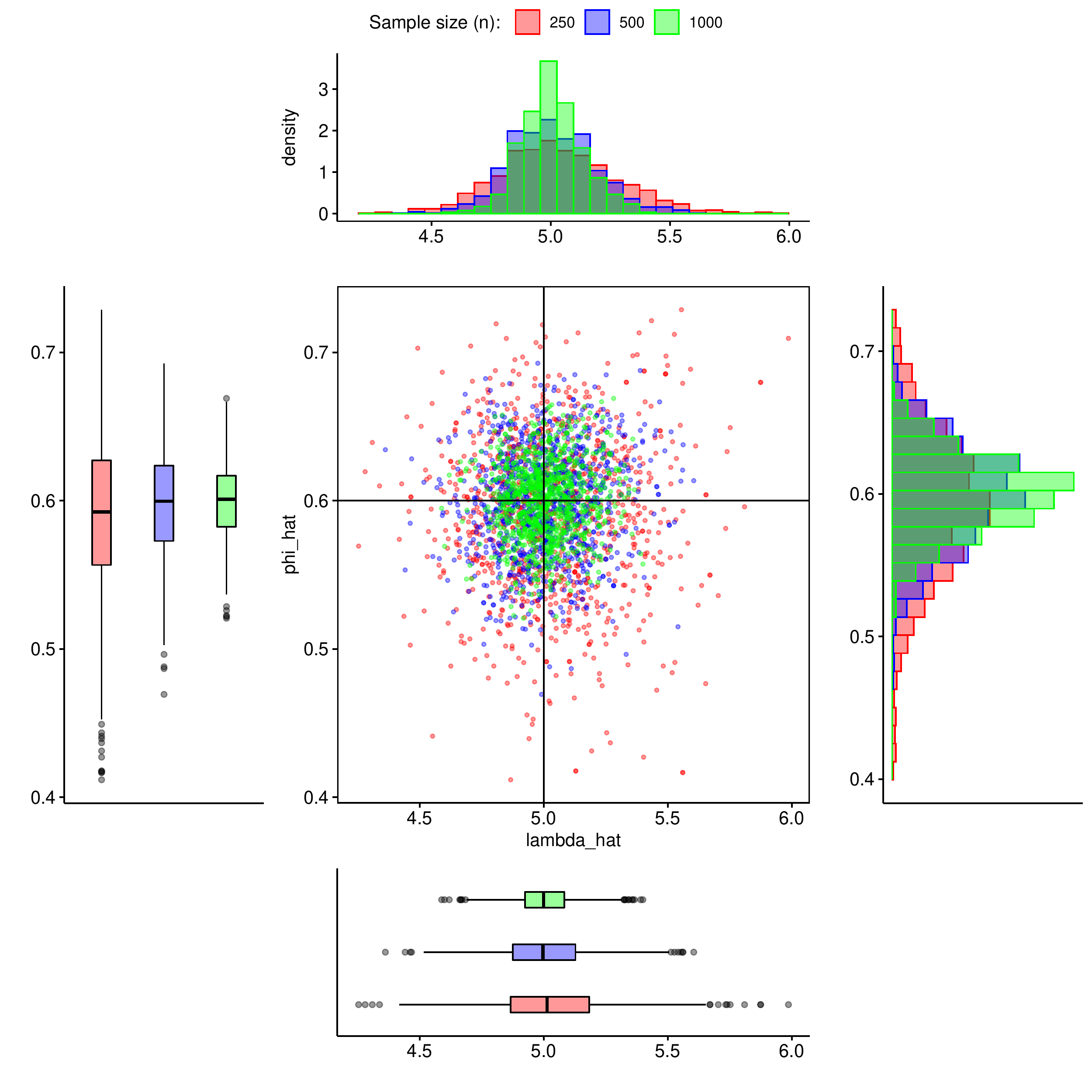}}
\caption{Simulation results for $\rho=0.5$, $\lambda=5$, $\phi=0.6$ and $\theta=0.4$. Presented are the scatter plot, and marginal densities and box plot for $n\in\{250,500,1000\}$ considering pairs $(\hat\phi,\hat\theta)$ (top-left), $(\hat\theta,\hat\lambda)$ (top-right) and $(\hat\phi,\hat\lambda)$ (bottom).} \label{ar2}
\end{figure}

From the simulation results we can also check the large sample behavior of the PMLE. In Figure \ref{ar2} we present pairwise scatter plots of $(\hat\lambda,\hat\phi,\hat\theta)$ for $\lambda = 5$, $\rho=0.5$, $\phi=0.6$ and $\theta=0.4$ along with marginal density and box plots, for all sample sizes. We observe that the scatter plots show a very bivariate normal-like behavior in all cases and, as the sample size increases, the estimated points concentrate around the DGP. From the scatter plots we observe that the correlation between $\hat\phi$ and $\hat\theta$ is relatively strong, while the correlation of $\hat\lambda$ with $\hat\phi$ and $\hat\theta$ is very small, especially when compared to the latter. Also notice that the marginal densities are very symmetrical around the DGP and the shape concentrates toward the true value as the sample size increases. All plots present evidence of the consistency and asymptotic normality of the PMLE in the context of UWARMA models.

\subsection*{Finite sample forecasting study}

\begin{figure}[h!]
\centering
\includegraphics[width=0.6\textwidth]{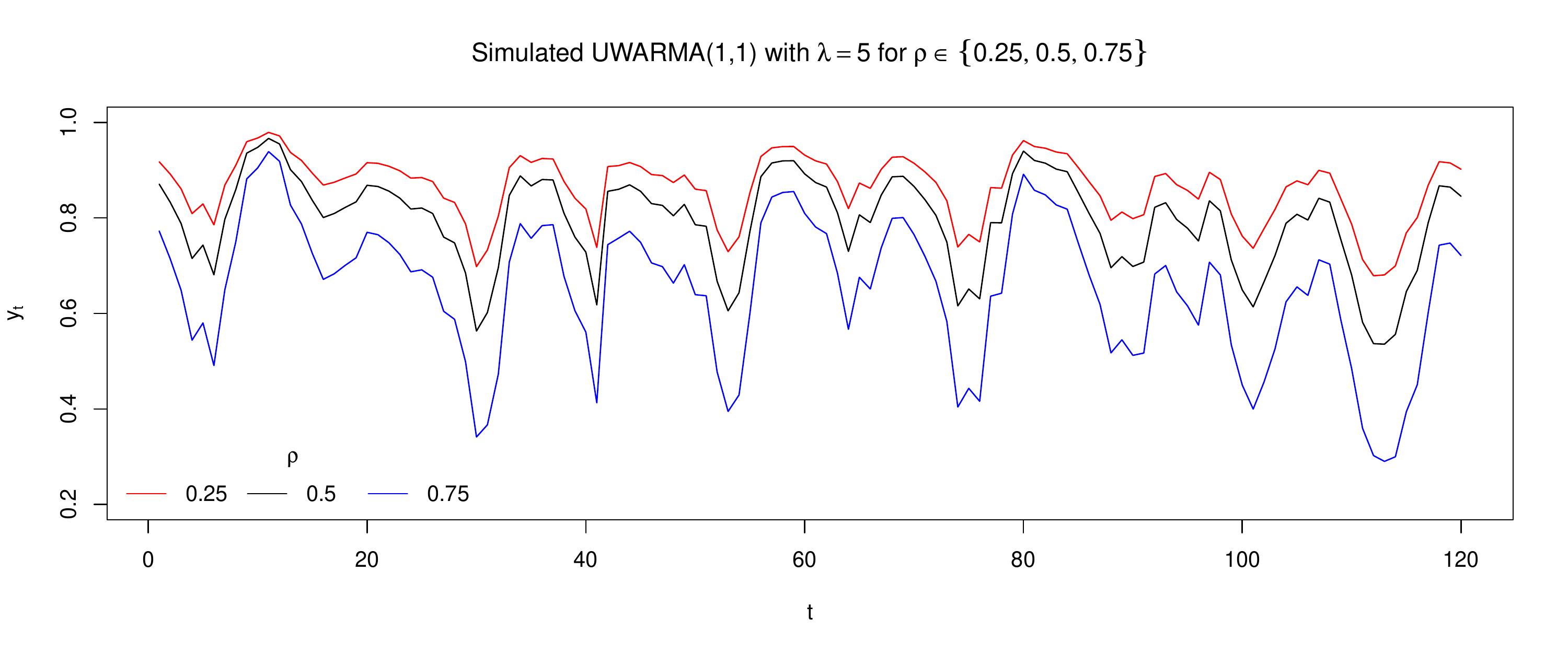}
\caption{Example of UWARMA$(1,1)$ with $\phi=0.6$, $\theta=0.4$, $\lambda = 5$ with covariates for different values of $\rho$ produced using the same random seed.} \label{cov}
\end{figure}

We now turn our attention to the forecasting capabilities of the proposed model. We simulate 1,000 replicas of an UWARMA$(1,1)$ with
$\phi=0.6$, $\theta=0.4$, $\lambda\in\{5,10,20\}$, $\rho\in\{0.25,0.5,0.75\}$ and $n=1{,}000$. We also include an intercept $\alpha=0.5$ and a single covariate $x_t:=\sin\Big(\frac{2\pi (t-6)}{12}\Big)$ with coefficient $\beta=0.5$. Figure \ref{cov} presents a typical sample path of the simulated model, for the different values of $\rho$, considering the same random seed.
In this exercise, for each generated time series we produce 24 step-ahead forecasts using \eqref{for} and calculate the mean absolute percentage error (MAPE) of $h$-step-ahead forecasts for $h\in\{1,6,12,18,24\}$. Table \ref{tfore} presents the average MAPE (in unitary value) along with the estimation of the parameter included in the model. From Table \ref{tfore} we notice that the MAPE are in general smaller than 10\% except for $\lambda=5$ and $\rho=0.75$. As expected the MAPE slowly increases with the horizon with few exceptions. Overall $\alpha$ and $\beta$ are always well estimated and so are $\phi$ and $\theta$, except for $\lambda = 20$, where the estimation of $\phi$ and $\theta$ present considerable bias, although the overall forecasting performance in these cases remain very good.
\begin{table}[h!]
\caption{Forecasting exercise: estimated values for the parameters $\alpha=\beta=0.5$, $\phi=0.6$, $\theta=0.4$ and $\lambda\in\{5,10,20\}$ and MAPE calculated for horizons $h\in\{1,6,12,18,24\}$.}\label{tfore}
\centering
\scriptsize
\vspace{.3cm}
\begin{tabular}{c|c|ccccr|ccccc}
%\hline
%\multicolumn{11}{c}{$\rho=0.25$}\\
\hline
$\rho$&$\lambda$& $\hat\alpha$  & $\hat\beta$ & $\hat\phi$ & $\hat\theta$ &\multicolumn{1}{c|}{$\hat\lambda$} & $h=1$ & $h=6$ & $h=12$& $h=18$ &$h=24$\\
\hline
\multirow{3}{*}{\begin{sideways}0.25\end{sideways}}
& 5  &  0.509  &  0.508  &  0.598  &  0.373  &  4.957  &  0.051  &  0.073  &  0.084  &  0.093  &  0.093    \\
& 10  &  0.490 &  0.506  &  0.610  &  0.319  &  9.743  &  0.033  &  0.046  &  0.052  &  0.058  &  0.058    \\
& 20  &  0.435  &  0.506  &  0.649  &  0.202  &  18.369  &  0.020  &  0.027  &  0.030  &  0.035  &  0.036    \\
%\hline
%\multicolumn{11}{c}{$\rho=0.5$}\\
\hline
 \multirow{3}{*}{\begin{sideways}0.5\end{sideways}}
 & 5  &  0.509  &  0.507  &  0.595  &  0.390  &  4.988  &  0.068  &  0.072  &  0.072  &  0.076  &  0.074    \\
 & 10  &  0.501  &  0.505  &  0.600  &  0.356  &  9.829  &  0.036  &  0.037  &  0.037  &  0.039  &  0.038    \\
 & 20  &  0.450  &  0.504  &  0.642  &  0.220  &  18.752  &  0.018  &  0.019  &  0.019  &  0.020  &  0.020    \\
%\hline
%\multicolumn{11}{c}{$\rho=0.75$}\\
\hline
\multirow{3}{*}{\begin{sideways}0.75\end{sideways}}
&  5  &  0.501  &  0.501  &  0.599  &  0.398  &  5.011  &  0.172  &  0.222  &  0.226  &  0.252  &  0.246   \\
&  10  &  0.508  &  0.504  &  0.593  &  0.390  &  9.937  &  0.063  &  0.079  &  0.081  &  0.087  &  0.085   \\
&  20  &  0.470   &  0.504  &  0.626  &  0.264  &  19.027  &  0.027  &  0.034  &  0.036  &  0.039  &  0.039 \\
\hline
\end{tabular}
\end{table}
\nopagebreak
\FloatBarrier
\section{Empirical Study}\label{empi}

Here, allied with monthly economic indicators, we evaluate our proposed model against KARMA and $\beta$ARMA models in an exercise of forecasting the Manufacturing Capacity Utilization, or CapU, from the US. According to \cite{Ragan1976}, the Federal Reserve combines surveys concerning the potential output and the actual production to estimate the CapU. Lying between 0 and 100 per cent, Capacity Utilization is defined as a ratio of the current output level to potential output. Because $\beta$ARMA does not allow for quantile analysis, and KARMA is only defined for the 50\% percentile, we focus only on the forecast of the median for our proposed method. The study of the effects of the covariates at each quantile will be left for future work.

It is not difficult to find prediction studies on the industry's output level, but only a handful use CapU as the predicted variable. For example, \cite{Baghestani2008} uses ARMA/ARMAX models, and \cite{TURHAN2015286} use Mixed data sampling (MIDAS) models to predict CapU. However, none of those methods guarantees that the prediction is bounded above and below. On the other hand, CapU is widely used as a predictor of inflation, but not the other way around. For example, \cite{CorradoMattey1997} show that inflation begins to accelerate whenever CapU exceeds the 82 percent mark. See \cite{bauer1990reexamination}, \cite{garner1994capacity} and \cite{ROSSI2010808}, for other studies relating inflation and CapU.

\begin{figure}[h!]
	\centering
\mbox{
	\includegraphics[width=0.5\textwidth]{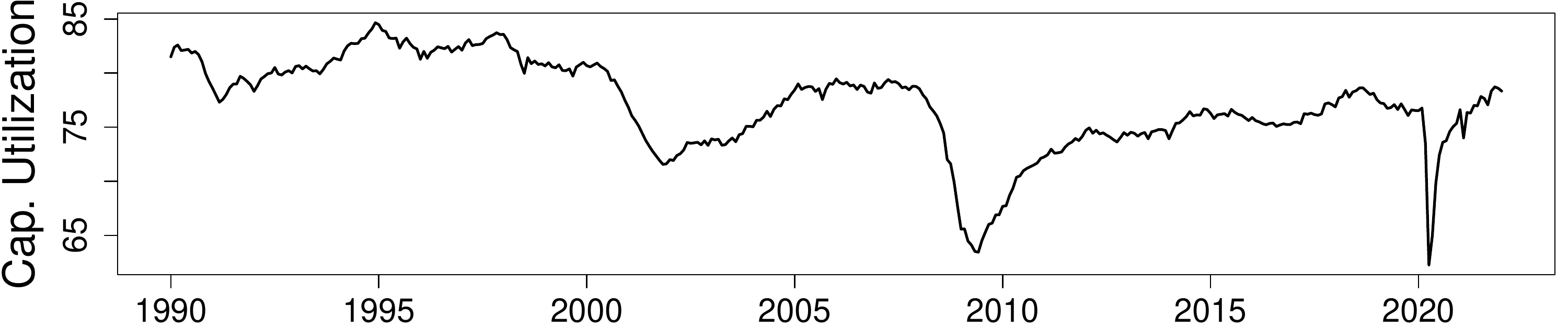}
	\includegraphics[width=0.5\textwidth]{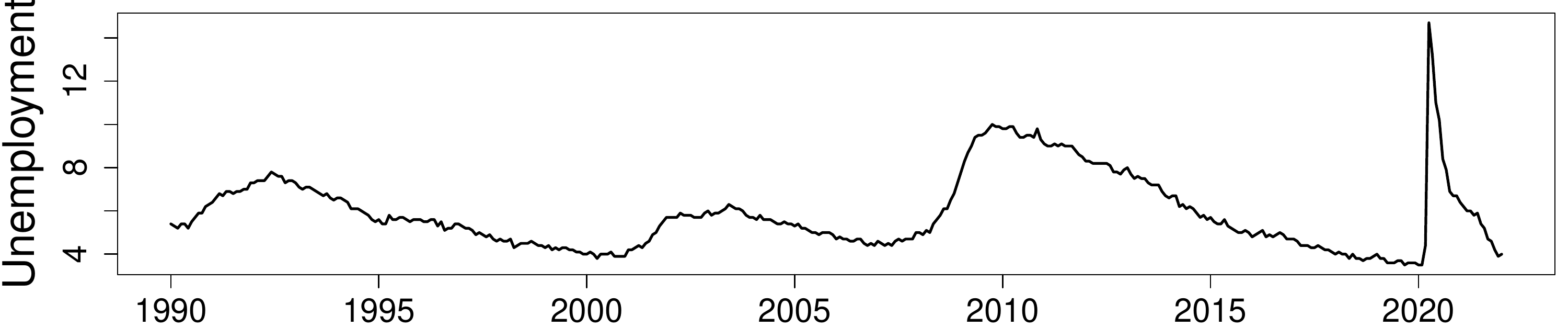}}
\mbox{
	\includegraphics[width=0.5\textwidth]{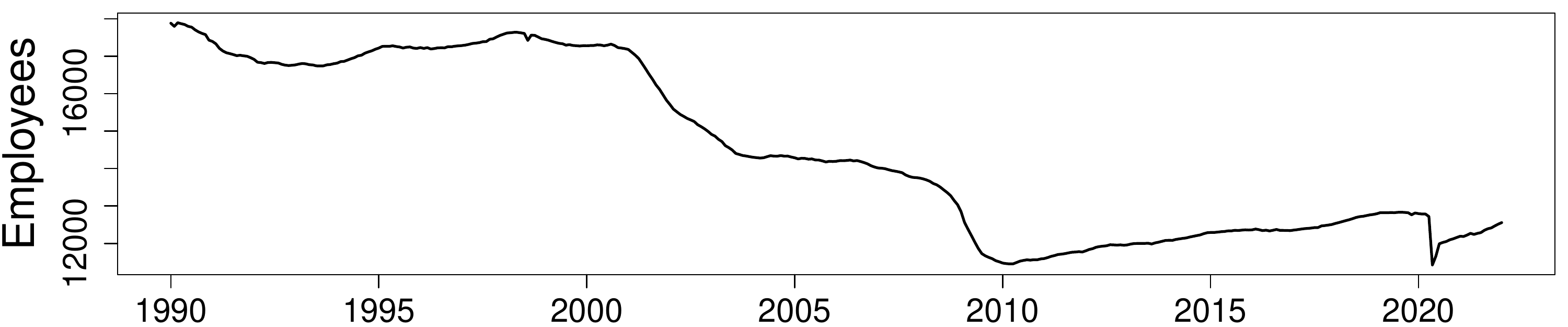}
	\includegraphics[width=0.5\textwidth]{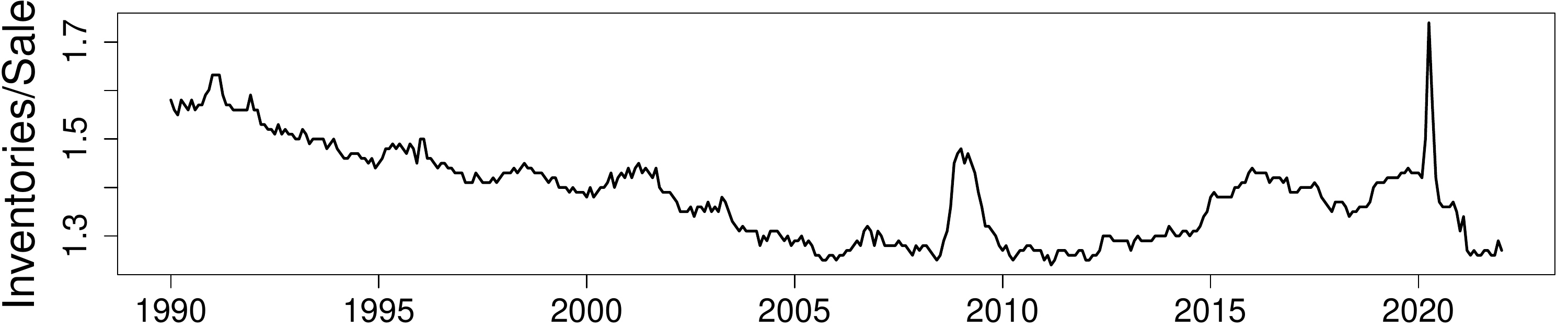}}
\mbox{
    \includegraphics[width=0.5\textwidth]{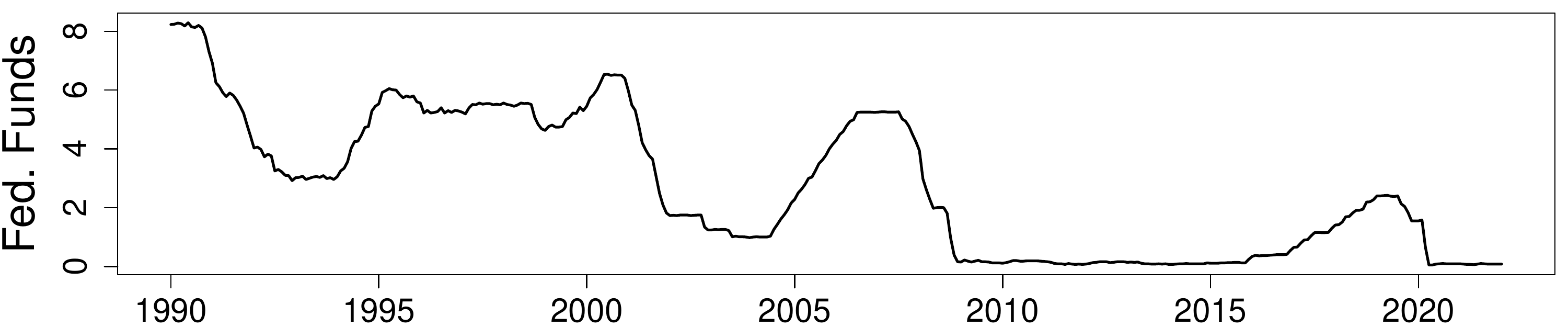}
	\includegraphics[width=   0.5 \textwidth]{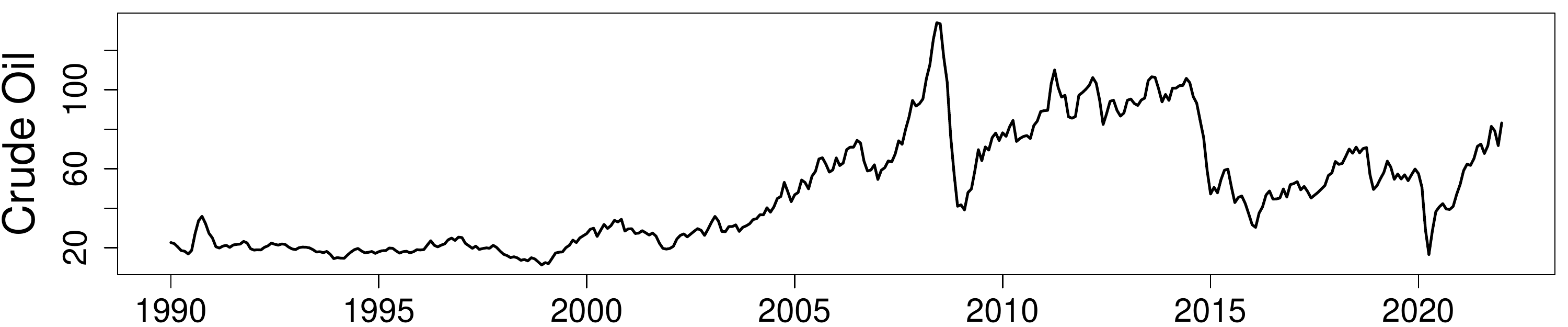}}
\mbox{
	\includegraphics[width=   0.5 \textwidth]{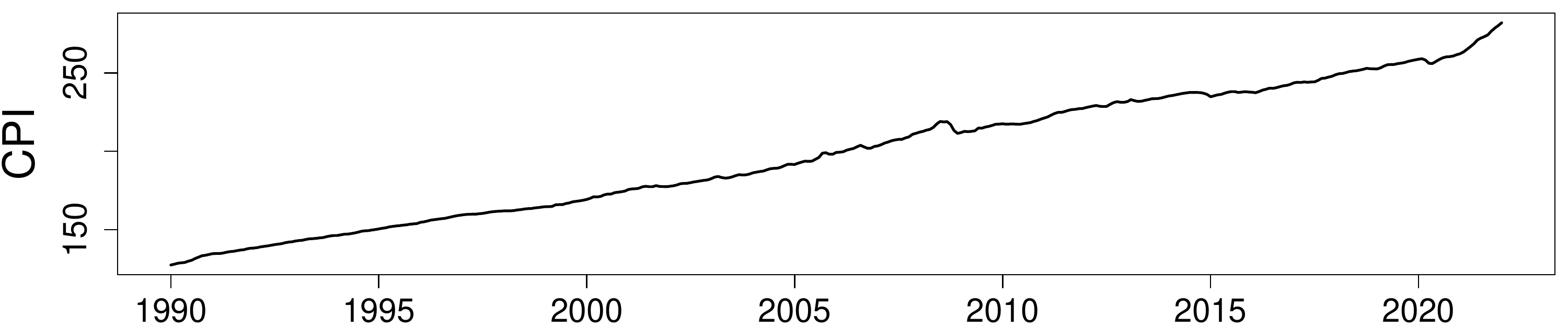}
	\includegraphics[width=   0.5 \textwidth]{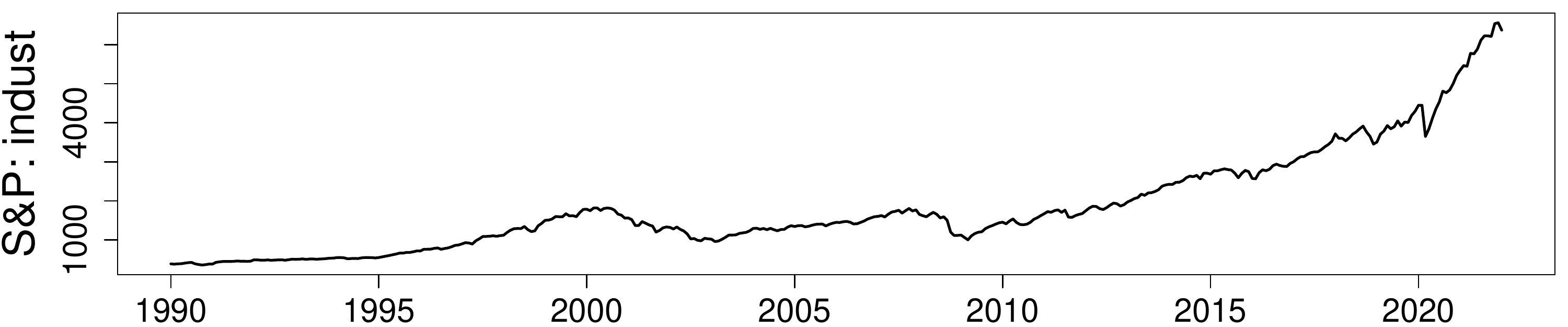}}
	\caption{Time series plot of Manufacturing Capacity Utilization (percent),
Civilian Unemployment Rate (percent), All Employees: Manufacturing (Thousands of Persons), and Total Business: Inventories to Sales Ratio (ratio), Effective Federal Funds Rate (percent), Crude Oil, spliced WTI and Cushing (dollars per barrel), CPI : All Items (index 1982-1984=100), and S\&P's Common Stock Price Index: Industrials (index).}
	\label{Empirical1}
\end{figure}

In this study, we use monthly data from January 1990 up to May 2022, totalizing 389 observations. The covariates used in this study are presented in Table \ref{tabcov}. Figure \ref{Empirical1} presents the evolution of each series through time. All series are public and available on the St. Louis FED website\footnote{https://fred.stlouisfed.org/}.

\begin{table}[htbp]
	\centering
	\caption{List and description of the exogenous variables applied in the study.}	
	\begin{tabular*}{1\textwidth}{@{\extracolsep{\fill} }cll}
		\toprule
		tcode & \multicolumn{1}{l}{fred} & \multicolumn{1}{l}{Description} \\
		\midrule
		2     & UNRATE & Civilian Unemployment Rate \\
		5     & MANEMP & All Employees: Manufacturing \\
		2     & ISRATIOx & Total Business: Inventories to Sales Ratio \\
		2     & FEDFUNDS & Effective Federal Funds Rate \\
		6     & OILPRICEx & Crude Oil, spliced WTI and Cushing \\
		6     & CPIAUCSL & CPI: All Items \\
		5     & S\&P: indust & S\&P's Common Stock Price Index: Industrials \\
		\bottomrule
	\end{tabular*}%
	\begin{tablenotes}[flushleft]
		\footnotesize
		\item \hspace{-.1cm}{Note:} tcode refers to the transformation necessary for the variable to become stationary, see \cite{McCracken2016}. Where (2) $\Delta x_t$, (5) $\Delta \log (x_t)$, and (6) $\Delta^2 \log (x_t)$.
	\end{tablenotes}
	\label{tabcov}%
\end{table}%

At any $t$, we use a rolling window of 287 months and forecast up to 6 months ahead, with 91 total months available for the forecast. The stationary versions of each covariates are included from lag 1 to 3. Given an ARMA($p,q$) structure, variable selection was performed using a backward elimination strategy based on $p$-values. In other words, after each estimation, the variable with the highest $p$-value greater than 0.05 is removed. This procedure is repeated until all included exogenous variables present a $p$-value lower than 0.05. We perform this procedure for all combinations of $p,q \in\{0,\cdots,3\}$, except for $p=q=0$. Nevertheless, we only report the configuration yielding the best results regarding MAPE for each type of model, namely, UWARMA$(2,0)$, KARMA$(1,0)$, and $\beta$ARMA$(3,2)$. Table \ref{tabresultsMAPE} shows that UWARMA(2,0) presents the lowest Average MAPE for every horizon. Conversely, KARMA(1,0) has the highest results for every horizon. This indicates that the UWARMA method is generally better than the other methods evaluated here.

%Nevertheless, we only report the configuration yielding the best results regarding RMSE and MAPE for each type of model, namely, UWARMA$(2,0)$, KARMA$(1,0)$, and $\beta$ARMA$(3,2)$.

%\begin{table}[htbp]
%	\centering
%	\caption{Average RMSE forecast results, multiplied by 100.}
%	\begin{tabular*}{1\textwidth}{@{\extracolsep{\fill} }lcccccc}
%		\hline
%		& \multicolumn{1}{c}{ $  t+1 $ } & \multicolumn{1}{c}{  $ t+2 $ } & \multicolumn{1}{c}{  $ t+3 $ } & \multicolumn{1}{c}{  $ t+4 $ } & \multicolumn{1}{c}{  $ t+5 $ } & \multicolumn{1}{c}{ $ t+6 $} \\
%		\hline
%		UWARMA(2,0)          & 0.0531 & 0.0694 & 0.0828 & 0.0882 & 0.0891 & 0.0924 \\
%		KARMA(1,0)           & 0.0600 & 0.0718 & 0.0773 & 0.0810 & 0.0826 & 0.0845 \\
%		$\beta$ARMA(3,2)     & 0.0589 & 0.0940 & 0.1160 & 0.1194 & 0.1233 & 0.1282 \\
%		\hline
%	\end{tabular*}%
%	%\begin{tablenotes}[flushleft]
%	%	\footnotesize
%	%	\item \hspace{-.1cm}{Note:} We only report, here, the best results for each method.
%	%\end{tablenotes}
%	\label{tabresultsRMSE}%
%\end{table}%

\begin{table}[htbp]
	\centering
	\caption{Average MAPE forecast results.}
	\begin{tabular*}{1\textwidth}{@{\extracolsep{\fill} }lcccccc}
		\toprule
		& \multicolumn{1}{c}{ $  t+1 $ } & \multicolumn{1}{c}{  $ t+2 $ } & \multicolumn{1}{c}{  $ t+3 $ } & \multicolumn{1}{c}{  $ t+4 $ } & \multicolumn{1}{c}{  $ t+5 $ } & \multicolumn{1}{c}{ $ t+6 $} \\
		\midrule
		UWARMA(2,0)          & 0.0135 & 0.0157 & 0.0181 & 0.0193 & 0.0209 & 0.0223 \\
		KARMA(1,0)           & 0.0201 & 0.0234 & 0.0251 & 0.0259 & 0.0267 & 0.0270 \\
		$\beta$ARMA(3,2)     & 0.0147 & 0.0182 & 0.0208 & 0.0216 & 0.0233 & 0.0251 \\
		\bottomrule
	\end{tabular*}%
	%	\begin{tablenotes}[flushleft]
		%		\footnotesize
		%		\item \hspace{-.1cm}{Note:} We only report, here, the best results for each method.
		%	\end{tablenotes}
	\label{tabresultsMAPE}%
\end{table}%

% Ajustar o testo com os novos resultados
%Tables \ref{tabresultsRMSE} and \ref{tabresultsMAPE} present the results in terms of RMSE and MAPE, respectively. Table \ref{tabresultsRMSE} shows a similar Average RMSE for UWARMA(2,0) and KARMA(1,0), but, except for $t+1$, $\beta$ARMA(3,2) produces a result reasonably higher than the other two methods. For the Average MAPE, $\beta$ARMA(3,2) shows similar outcomes to the other methods, but now UWARMA(2,0) presents the lowest values for every horizon. The results indicate that the UWARMA method is equal, and in some cases better, than the other methods here evaluated.

Concerning the use of variables and their lags, UWARMA(2,0) is the most parsimonious, using, on average, 4.7 covariates from the available ones. KARMA(1,0) and $\beta$ARMA(3,2) use, on average, 7.5 and 8.2, respectively. Table \ref{tabfreq} shows the frequency selection rate of each variable, i.e., the average of whether at least one lag is chosen. S\&P: indust is always significant for $\beta$ARMA(3,2) and also de the most important variable for UWARMA(2,0), selected 94\% of the time. For KARMA(1,0), MANEMP is the most important variable, which is chosen 74\% of the time. Except for S\&P: indust, given the different variable selection rates between models, there is no consensus on the order of the most influential variables across models.

\begin{table}[htbp]
	\centering
	\caption{Frequency selection rate of each variable, given the model.}
	\begin{tabular*}{1\textwidth}{@{\extracolsep{\fill} }lccc}
		\toprule
		& \multicolumn{1}{l}{UWARMA(2,0)} & \multicolumn{1}{l}{KARMA(1,0)} & \multicolumn{1}{l}{$\beta$ARMA(3,2)} \\
		\midrule
		UNRATE & 0.8681 & 0.6813 & 0.9341 \\
		MANEMP & 0.7143 & 0.7363 & 0.0110 \\
		ISRATIOx & 0.2527 & 0.5604 & 0.9121 \\
		FEDFUNDS & 0.1648 & 0.7033 & 0.9780 \\
		OILPRICEx & 0.0659 & 0.5714 & 0.4945 \\
		CPIAUCSL & 0.8242 & 0.3626 & 0.8571 \\
		S\&P: indust & 0.9451 & 0.6923 & 1.0000 \\
		\bottomrule
	\end{tabular*}%
	\label{tabfreq}%
\end{table}%
\FloatBarrier
\bibliographystyle{apalike}
\bibliography{UW}

\begin{thebibliography}{}

\bibitem[Baghestani, 2008]{Baghestani2008}
Baghestani, H. (2008).
\newblock Predicting capacity utilization: Federal reserve vs time-series
  models.
\newblock {\em J. Econ. Fin.}, 32(1):47--57.

\bibitem[Bauer et~al., 1990]{bauer1990reexamination}
Bauer, P.~W. et~al. (1990).
\newblock A reexamination of the relationship between capacity utilization and
  inflation.
\newblock {\em Economic Review}, 26(2):2--12.

\bibitem[Bayer et~al., 2017]{Bayers}
Bayer, F.~M., Bayer, D.~M., and Pumi, G. (2017).
\newblock Kumaraswamy autoregressive moving average models for double bounded
  environmental data.
\newblock {\em Journal of Hydrology}, 555:385--396.

\bibitem[Bayer et~al., 2018]{Bayerseas}
Bayer, F.~M., Cintra, R.~J., and Cribari-Neto, F. (2018).
\newblock Beta seasonal autoregressive moving average models.
\newblock {\em Journal of Statistical Computation and Simulation},
  88(15):2961--2981.

\bibitem[Benaduce and Pumi, 2022]{helen}
Benaduce, H.~S. and Pumi, G. (2022).
\newblock {SYMARFIMA}: a dynamical model for conditionally symmetric time
  series with long range dependence mean structure.
\newblock {\em Journal of Statistical Planning and Inference}, Forthcoming.

\bibitem[Benjamin et~al., 2003]{benjamin2003}
Benjamin, M.~A., Rigby, R.~A., and Stasinopoulos, D.~M. (2003).
\newblock Generalized autoregressive moving average models.
\newblock {\em Journal of the American Statistical Association},
  98(461):214--223.

\bibitem[Casarin et~al., 2012]{Casarin2012}
Casarin, R., Dalla~Valle, L., and Leisen, F. (2012).
\newblock Bayesian model selection for beta autoregressive processes.
\newblock {\em Bayesian Analysis}, 7(2):385--410.

\bibitem[Corrado and Mattey, 1997]{CorradoMattey1997}
Corrado, C. and Mattey, J. (1997).
\newblock Capacity utilization.
\newblock {\em J. Econ. Perspect.}, 11(1):151--167.

\bibitem[Cox et~al., 1981]{cox1981}
Cox, D.~R., Gudmundsson, G., Lindgren, G., Bondesson, L., Harsaae, E., Laake,
  P., Juselius, K., and Lauritzen, S.~L. (1981).
\newblock Statistical analysis of time series: Some recent developments [with
  discussion and reply].
\newblock {\em Scandinavian Journal of Statistics}, pages 93--115.

\bibitem[Ferrari and Cribari-Neto, 2004]{Ferrari2004}
Ferrari, S. L.~P. and Cribari-Neto, F. (2004).
\newblock Beta regression for modelling rates and proportions.
\newblock {\em Journal of Applied Statistics}, 31(7):799--815.

\bibitem[Fokianos and Kedem, 1998]{Fokianos1998}
Fokianos, K. and Kedem, B. (1998).
\newblock Prediction and classification of non-stationary categorical time
  series.
\newblock {\em Journal of Multivariate Analysis}, 67:277--296.

\bibitem[Fokianos and Kedem, 2004]{Fokianos2004}
Fokianos, K. and Kedem, B. (2004).
\newblock Partial likelihood inference for time series following generalized
  linear models.
\newblock {\em Journal of Time Series Analysis}, 25(2):173--197.

\bibitem[Garner et~al., 1994]{garner1994capacity}
Garner, C.~A. et~al. (1994).
\newblock Capacity utilization and us inflation.
\newblock {\em Economic Review-Federal Reserve Bank of Kansas City}, 79:5--5.

\bibitem[Gradshteyn and Ryzhik, 2007]{GR2007}
Gradshteyn, I.~S. and Ryzhik, I.~M. (2007).
\newblock {\em Table of integrals, series, and products}.
\newblock Academic Press, 7 edition.

\bibitem[Grande et~al., 2022]{grande}
Grande, A.~F., Pumi, G., and Cybis, G.~B. (2022).
\newblock Granger causality and time series regression for modeling the
  migratory dynamics of {I}nfluenza into {B}razil.
\newblock {\em SORT}, Forthcoming.

\bibitem[Kedem and Fokianos, 2002]{Kedem2002}
Kedem, B. and Fokianos, K. (2002).
\newblock {\em Regression models for time series analysis}.
\newblock John Wiley \& Sons.

\bibitem[Maior and Cysneiros, 2018]{Maior}
Maior, V. and Cysneiros, F. (2018).
\newblock {SYMARMA}: a new dynamic model for temporal data on conditional
  symmetric distribution.
\newblock {\em Statistical Papers}, 59.

\bibitem[Mazucheli et~al., 2020]{UWreg}
Mazucheli, J., Menezes, A. F.~B., Fernandes, L.~B., de~Oliveira, R.~P., and
  Ghitany, M.~E. (2020).
\newblock The unit-{W}eibull distribution as an alternative to the
  {K}umaraswamy distribution for the modeling of quantiles conditional on
  covariates.
\newblock {\em Journal of Applied Statistics}, 47(6):954--974.

\bibitem[McCracken and Ng, 2016]{McCracken2016}
McCracken, M.~W. and Ng, S. (2016).
\newblock {FRED-MD}: A monthly database for macroeconomic research.
\newblock {\em Journal of Business and Economic Statistics}, 34(4):574--589.

\bibitem[Mitnik and Baek, 2013]{MK}
Mitnik, P.~A. and Baek, S. (2013).
\newblock The {K}umaraswamy distribution: median-dispersion
  re-parameterizations for regression modeling and simulation-based estimation.
\newblock {\em Statistical Papers}, 54(1):177--192.

\bibitem[Prass et~al., 2022]{Prass}
Prass, T.~S., Carlos, J.~H., Taufemback, C.~G., and Pumi, G. (2022).
\newblock Positive time series regression models.

\bibitem[Pumi et~al., 2021]{BARC}
Pumi, G., Prass, T.~S., and Souza, R.~R. (2021).
\newblock A dynamic model for double-bounded time series with chaotic-driven
  conditional averages.
\newblock {\em Scandinavian Journal of Statistics}, 48(1):68--86.

\bibitem[Pumi et~al., 2020]{Pumi2020}
Pumi, G., Rauber, C., and Bayer, F.~M. (2020).
\newblock Kumaraswamy regression model with {A}randa-{O}rdaz link function.
\newblock {\em TEST}, 29:1051--1071.

\bibitem[Pumi et~al., 2019]{Pumi2017}
Pumi, G., Valk, M., Bisognin, C., Bayer, F.~M., and Prass, T.~S. (2019).
\newblock Beta autoregressive fractionally integrated moving average models.
\newblock {\em Journal of Statistical Planning and Inference}, 200:196--212.

\bibitem[{R Core Team}, 2022]{R}
{R Core Team} (2022).
\newblock {\em R: A Language and Environment for Statistical Computing}.
\newblock R Foundation for Statistical Computing, Vienna, Austria.

\bibitem[Ragan, 1976]{Ragan1976}
Ragan, J.~F. (1976).
\newblock Measuring capacity utilization in manufacturing.
\newblock {\em Federal Reserve Board New York Quarterly Review}, 1:13--28.

\bibitem[Rocha and Cribari-Neto, 2009]{Rocha2009}
Rocha, A.~V. and Cribari-Neto, F. (2009).
\newblock Beta autoregressive moving average models.
\newblock {\em Test}, 18(3):529--545.

\bibitem[Rocha, 2017]{Rocha2017}
Rocha, A. V.and Cribari-Neto, F. (2017).
\newblock Erratum to: Beta autoregressive moving average models.
\newblock {\em Test}, 26(2):451--459.

\bibitem[Rossi and Sekhposyan, 2010]{ROSSI2010808}
Rossi, B. and Sekhposyan, T. (2010).
\newblock Have economic models' forecasting performance for {US} output growth
  and inflation changed over time, and when?
\newblock {\em International Journal of Forecasting}, 26(4):808--835.

\bibitem[Scher et~al., 2020]{port}
Scher, V.~T., Cribari-Neto, F., Pumi, G., and Bayer, F.~M. (2020).
\newblock Goodness-of-fit tests for $\beta${ARMA} hydrological time series
  modeling.
\newblock {\em Environmetrics}, 31(3):e2607.

\bibitem[Turhan et~al., 2015]{TURHAN2015286}
Turhan, I.~M., Sensoy, A., and Hacihasanoglu, E. (2015).
\newblock Shaping the manufacturing industry performance: Midas approach.
\newblock {\em Chaos, Solitons \& Fractals}, 77:286--290.

\end{thebibliography}

\section*{Appendix}

\begin{lema}\label{lema}
Let $Y\sim \mathrm{UW}(\mu,\lambda;\rho)$ for $\rho,\mu\in(0,1)$ and $\lambda>0$. Then
\begin{align}
&\E\bigg(\bigg[\frac{\log(Y)}{\log(\mu)}\bigg]^{\lambda}\bigg)=-\frac{1}{\log(\rho)},\label{L1} \\
&\E\bigg(\log\bigg(\frac{\log(Y)}{\log(\mu)}\bigg)\bigg)=\frac{\kappa+\log\big(-\log(\rho)\big)}{\lambda},\label{L2}\\
&\E\bigg(\bigg[\frac{\log(Y)}{\log(\mu)}\bigg]^{\lambda}\log\bigg(\frac{\log(Y)}{\log(\mu)}\bigg)\bigg)=
\frac{\kappa+\log\big(-\log(\rho)\big)-1}{\lambda\log(\rho)},\label{L3}\\
&\E\bigg(\bigg[\frac{\log(Y)}{\log(\mu)}\bigg]^{\lambda}\log\bigg(\frac{\log(Y)}{\log(\mu)}\bigg)^2\bigg)=
\frac{\pi^2+6(\kappa-2)\kappa-6\log\bigl(-\log(\rho)\bigr)\bigl[\log\bigl(-\log(\rho)\bigr)+2\kappa-2\bigr]}{6\lambda^2\log(\rho)}, \label{L4}
\end{align}
where $\kappa=0.5772156649\dots$  is the Euler-Mascheroni constant \citep{GR2007}.
\end{lema}
\noindent\textbf{Proof:} Conveniently writing
\begin{align*}
\E\bigg(\bigg[\frac{\log(Y)}{\log(\mu)}\bigg]^{\lambda}\bigg) = \log(\rho)\int_0^1 \bigg[\frac{\log(y)}{\log(\mu)}\bigg]^{\lambda}\rho^{{\big(\frac{\log(y)}{\log(\mu)}\big)}^\lambda}
\frac{\lambda}{y\log(y)}\bigg[\frac{\log(y)}{\log(\mu)}\bigg]^{\lambda}dy,
\end{align*}
and changing variables to
\begin{equation}\label{cv}
u=\bigg[\frac{\log(y)}{\log(\mu)}\bigg]^{\lambda} \quad \Longrightarrow\quad du=\frac{\lambda}{y\log(y)}\bigg[\frac{\log(y)}{\log(\mu)}\bigg]^{\lambda}dy,
\end{equation}
we obtain
\begin{equation*}\E\bigg(\bigg[\frac{\log(Y)}{\log(\mu)}\bigg]^{\lambda}\bigg) = -\log(\rho)\int_0^\infty u\rho^udu=-\log(\rho)\bigg[\frac{(u\log(\rho)-1)\rho^u}{\log(\rho)^2}\bigg|_{u=0}^\infty\bigg]=-\frac{1}{\log(\rho)},\end{equation*}
where the last equality follows since $0<\rho<1$. This proves \eqref{L1} and we move to prove \eqref{L2}. Applying the change of variables \eqref{cv}, we have
\begin{equation*}\E\bigg(\log\bigg(\frac{\log(Y)}{\log(\mu)}\bigg)\bigg)=\frac{\log(\rho)}{\lambda}\int_0^\infty\log(z)\rho^zdz.\end{equation*}
By using formula 2.751.2 in \cite{GR2007}, we have
\begin{align}\label{int}
\int_0^\infty \log(z)\rho^zdz &= \frac1{\log(\rho)}\bigg(\rho^z\log(z)-\mathrm{Ei}\big(\log(\rho)z\big)\bigg|_{z=0}^\infty\bigg)\nonumber\\
&= \frac1{\log(\rho)} \Big(-\lim_{z\rightarrow 0^+}\big[\rho^z\log(z)-\mathrm{Ei}\big(\log(\rho)z\big)\big]\Big)\nonumber\\
&=\frac{\kappa+\log\big(-\log(\rho)\big)}{\log(\rho)},
\end{align}
where Ei$(x):=-\int_{-x}^\infty \frac{e^{-t}}{t}dt$ is the exponential integral function. The limit above is obtained by using formula 8.214.1 in \cite{GR2007} to write
\begin{equation*}
\rho^z\log(z)-\mathrm{Ei}\big(\log(\rho)z\big) = (\rho^z-1)\log(z)-\kappa-\log\big(-\log(\rho)\big)-\sum_{k=1}^\infty \frac{z^k\log(\rho)^k}{k\cdot k!},
\end{equation*}
which yields the desired limit and \eqref{L2} follows. To show \eqref{L3},  upon applying the change of variables \eqref{cv} once again, we obtain
\begin{equation}\label{esp}
\E\bigg(\bigg[\frac{\log(Y)}{\log(\mu)}\bigg]^{\lambda}\log\bigg(\frac{\log(Y)}{\log(\mu)}\bigg)\bigg) = -\frac{\log(\rho)}\lambda \int_0^\infty z\rho^z\log(z)dz.
\end{equation}
Integration by parts with
\begin{equation*}u=z\log(z)\ \Rightarrow\ du=(\log(z)+1)dz\qquad\mbox{and}\qquad dv=\rho^zdz\ \Rightarrow\ v=\frac{\rho^z}{\log(\rho)},\end{equation*}
yields
\begin{align*}
\int_0^\infty z\rho^z\log(z)dz &= \frac{z\log(z)\rho^z}{\log(\rho)}\bigg|_{z=0}^\infty - \frac1{\log(\rho)}\bigg[\int_0^\infty \log(z)\rho^zdz+\int_0^\infty\rho^zdz\bigg]\\
&=- \frac1{\log(\rho)}\bigg[\frac{\kappa+\log\big(-\log(\rho)\big)}{\log(\rho)}-\frac1{\log(\rho)}\bigg],
\end{align*}
by \eqref{int}. Now \eqref{L3} follows from \eqref{esp} and the result above. Finally, to show \eqref{L4}, we apply \eqref{cv} again to obtain
\begin{equation*}
\E\bigg(\bigg[\frac{\log(Y)}{\log(\mu)}\bigg]^{\lambda}\log\bigg(\frac{\log(Y)}{\log(\mu)}\bigg)^2\bigg)=-\frac{\log(\rho)}{\lambda^2}\int_0^\infty z\rho^z \log(z)^2dz.
\end{equation*}
integration by parts with
\begin{equation*}
u=z\log(z)^2 \ \Rightarrow \ du=\big[\log(z)^2+2\log(z)\big]dz\qquad\mbox{and}\qquad dv=\rho^zdz\ \Rightarrow\ v=\frac{\rho^z}{\log(\rho)},
\end{equation*}
yields
\begin{equation}\label{fin}
\int_0^\infty z\rho^z \log(z)^2dz = \frac{z\log(z)^2\rho^z}{\log(\rho)}\bigg|_{z=0}^\infty - \frac1{\log(\rho)}\int_{0}^\infty\rho^{x}\log(x)^2dx-
\frac2{\log(\rho)}\int_{0}^\infty\rho^{x}\log(x)dx.
\end{equation}
The first term on the right hand side of \eqref{fin} is 0, while the third one has been computed in \eqref{int}. The second term can be evaluated similarly, yielding
\[\int_{0}^\infty\rho^{x}\log(x)^2dx=-\frac{6\log\big(-\log(\rho)\big)\big[\log\big(-\log(\rho)\big)+2\kappa\big]+\pi^2+6\kappa^2}{6\log(\rho)},\]
and \eqref{L4} follows. This completes the proof. \qed

\end{document}